%% file: main.tex
\documentclass[a4paper,review]{cas-sc}

\usepackage[section]{placeins} 
\input{preamble}

\begin{document}

\title[mode=title]{Utilizing Time-Reversibility for Shock Capturing in Nonlinear Hyperbolic Conservation Laws}
\shorttitle{Utilizing Time-Reversibility for Shock Capturing in Nonlinear Hyperbolic Conservation Laws}
\shortauthors{T. Dzanic \textit{et al.}}

\author[1]{T. Dzanic}[orcid=0000-0003-3791-1134]
\cormark[1]
\cortext[cor1]{Corresponding author}
\ead{tdzanic@tamu.edu}
\author[1,2]{W. Trojak}[orcid=0000-0002-4407-8956]
\author[1]{F. D. Witherden}[orcid=0000-0003-2343-412X]

\address[1]{Department of Ocean Engineering, Texas A\&M University, College Station, TX 77843}
\address[1]{Department of Aerospace Engineering, Texas A\&M University, College Station, TX 77843}

\begin{abstract}
In this work, we introduce a novel approach to formulating an artificial viscosity for shock capturing in nonlinear hyperbolic systems by utilizing the property that the solutions of hyperbolic conservation laws are not reversible in time in the vicinity of shocks. The proposed approach does not require any additional governing equations or \textit{a priori} knowledge of the hyperbolic system in question, is independent of the mesh and approximation order, and requires the use of only one tunable parameter. The primary novelty is that the resulting artificial viscosity is unique for each component of the conservation law which is advantageous for systems in which some components exhibit discontinuities while others do not. The efficacy of the method is shown in numerical experiments of multi-dimensional hyperbolic conservation laws such as nonlinear transport, Euler equations, and ideal magnetohydrodynamics using a high-order discontinuous spectral element method on unstructured grids.  

\end{abstract}

\begin{keywords}
Shock capturing \sep Artificial viscosity \sep Hyperbolic systems
\sep Euler equations \sep Ideal magnetohydrodynamics \sep Spectral element method
\end{keywords}



\maketitle

\input{introduction}

\input{methodology}

\input{analysis}
\input{implementation}
\input{results}

\input{conclusions}

\section*{Acknowledgements}
\label{sec:ack}
This work was supported in part by the U.S. Air Force Office of Scientific Research via grant FA9550-21-1-0190 (``Enabling next-generation heterogeneous computing for massively parallel high-order compressible CFD'') of the Defense University Research Instrumentation Program (DURIP) under the direction of Dr. Fariba Fahroo.

\bibliographystyle{unsrtnat}
\bibliography{reference}




\end{document}

%% file: preamble.tex
\usepackage{amssymb}
\usepackage{amsthm}
\usepackage{amsmath}
\usepackage{upgreek}
\usepackage{pdflscape}
\usepackage{listings}
\usepackage{multirow}
\usepackage[percent]{overpic}
\usepackage{multicol}
\usepackage{color}
\usepackage{stmaryrd}
\usepackage{xfrac}
\usepackage[capitalise]{cleveref}
\usepackage{booktabs}
\usepackage[section]{placeins} 
\usepackage[section]{algorithm}
\usepackage{calc}
\usepackage{mathtools}
\usepackage[titletoc]{appendix}
\usepackage{siunitx}
\usepackage[sort&compress,square,numbers]{natbib}

\usepackage{graphicx}
\usepackage{graphics}
\usepackage{wrapfig}
\usepackage{float}
\usepackage{subfig}
\usepackage[percent]{overpic}
\usepackage{varwidth}
\usepackage{tikz}
\usepackage{pgfplots}
\usepackage{adjustbox}
\usetikzlibrary{arrows,matrix,positioning,fit}
\usetikzlibrary{shapes,positioning}
\usetikzlibrary{backgrounds}
\usepackage{tikz-layers}
\usepgfplotslibrary{fillbetween}
\usetikzlibrary{intersections}
\pgfplotsset{compat=1.14}
\usepgfplotslibrary{colorbrewer}
\usepgfplotslibrary{patchplots}
\usepgfplotslibrary[colorbrewer]
\usetikzlibrary{pgfplots.colorbrewer}
\usetikzlibrary[pgfplots.colorbrewer]
\usepgfplotslibrary{units}
\usetikzlibrary{spy}
\usepackage{pgfplotstable}
\usepackage{arrayjobx}
\graphicspath{ {./figs/} }
\usetikzlibrary{external}

\newlength\myheight
\newlength\mydepth
\settototalheight\myheight{Xygp}
\newcommand*\inlinegraphics[1]{%
  \settototalheight\myheight{Xygp}%
  \settodepth\mydepth{Xygp}%
  \raisebox{-\mydepth}{\includegraphics[height=\myheight]{#1}}%
}
\newcommand\orcid[1]{\href{https://orcid.org/#1}{\inlinegraphics{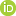}}}

\makeatletter
\def\BState{\State\hskip-\ALG@thistlm}
\makeatother

\newdefinition{definition}{Definition}[section]

\newcommand{\half}{{\frac{1}{2}}}
\newcommand{\shalf}{{\sfrac{1}{2}}}

%% file: introduction.tex
\section{Introduction}
\label{sec:intro}
The approximation of systems of nonlinear hyperbolic conservation laws remains a challenge for the development of many high-resolution numerical schemes due to their lack of robustness in the vicinity of discontinuities. For these systems, discontinuities can, and often do, arise from initially-smooth data \citep{Hopf1950}, and without proper treatment, can result in spurious oscillations that can degrade the quality of the solution or cause the scheme to fail altogether. Various nonlinear stabilization techniques have emerged among practitioners of high-order numerical schemes, and these techniques share a common goal of robustly resolving discontinuous features while recovering high-order accuracy in smooth regions. These approaches generally fall into one of four categories: limiting \citep{Guermond2019,Boris1997}, filtering \citep{Kanevsky2006,Glaubitz2017,Hesthaven2008}, stencil modification \citep{Jiang1996, Liu1994}, or artificial viscosity \citep{Persson2006,Guermond2011,Barter2007}, and each of these approaches differ in their accuracy, robustness, computational cost, and generalizability in the context of arbitrary unstructured grids. 

The addition of an artificial viscosity can be considered to be the most ubiquitous stabilization method for numerical schemes, dating as far back as to the work of \citet{VonNeumann1950}. This approach possesses several favorable properties in that it can robustly resolve discontinuities, is simple to implement and extend to unstructured grids, and has rigorous mathematical backing in the context of hyperbolic conservation laws \citep{Dafermos2010_4}. However, determining when and how much artificial dissipation to apply is an open problem, as the stability of the scheme can suffer if not enough dissipation is introduced in the vicinity of a shock, but solution accuracy in smooth regions can be severely degraded if unnecessary dissipation is applied. Previous works have taken both physical and numerical approaches for estimating this distribution of artificial viscosity. The use of the residual of the conservation law to form an artificial viscosity was explored from early works such as that of \citet{Johnson1990} to more recent works such as that of \citet{Nazarov2012}. This idea was extended to a more physically-relevant approach in the work of \citet{Guermond2011} where the violation of an entropy conservation equation (i.e., entropy residual) was used as an indicator for when to apply artificial dissipation, with the magnitude of the artificial viscosity proportional to the relative entropy production. From a purely numerical approximation approach, \citet{Persson2006} showed that the suboptimal convergence of the modal expansion of the solution in the vicinity of a shock can be used as an indicator for where to apply artificial dissipation. Approximate scaling arguments for the magnitude of the artificial viscosity were then used to quantify the amount of artificial dissipation to introduce. This approach was then used in \citet{Barter2007} to form a PDE-based artificial viscosity term via a forcing term driven by the shock indicator. These methods have shown varying levels of success across many applications, effectively differing in their accuracy, sensitivity to tunable parameters, and computational cost. 

In this work, we describe a novel approach for estimating the artificial dissipation required for robustly resolving discontinuous features in hyperbolic conservation laws. The artificial viscosity term is calculated based on an estimation of the \textit{reversibility} of the system, utilizing the fact that entropy production by shocks imposes a time-irreversibility condition on the solution \citep{Dafermos2010_4}. The proposed approach can be applied without any additional governing equations, is independent of the mesh and approximation order, does not require \textit{a priori} knowledge of the hyperbolic system in question, and requires the use of only one tunable parameter. The primary novelty of this method is that the resulting artificial viscosity is unique for each component of the conservation law which can significantly improve the accuracy of the approximation of systems in which some components exhibit discontinuities while others do not without additional computational cost. The efficacy of this method is shown through numerical simulations of hyperbolic conservation laws in the context of a high-order discontinuous spectral element method. 

The remainder of this paper is organized as follows. The proposed time-reversibility approach for artificial viscosity is introduced in \cref{sec:methodology}. The underlying mechanisms of this approach are then analyzed in \cref{sec:analysis}, and numerical implementation details for a high-order discontinuous spectral element method are given in \cref{sec:implementation}. The results of numerical experiments for a variety of multi-dimensional nonlinear hyperbolic conservation laws such as nonlinear transport, Euler equations, and ideal magnetohydrodynamics are shown in \cref{sec:results}, and conclusions are then drawn in \cref{sec:conclusion}.

%% file: methodology.tex
\section{Methodology}\label{sec:methodology}
Consider an arbitrary hyperbolic conservation law of the form
\begin{equation}\label{eq:hyp}
    \partial_t \mathbf{u} (\mathbf{x}, t) + \boldsymbol{\nabla}{\cdot} \mathbf{F}\left(\mathbf{u}(\mathbf{x}, t)\right) = 0, \quad \mathrm{for}\ (\mathbf{x}, t) \in \Omega \times \mathbb R^+,
\end{equation}
where $\mathbf{u} \in \mathbb R^m$ is a vector-valued solution, $m$ is an arbitrary number of field variables, $\mathbf{F}(\mathbf{u}) \in \mathbb R^{m\times d}$ is the flux, $\Omega \subset \mathbb R^d$ is the domain, and $d$ is some arbitrary spatial dimension. For simplicity, the domain is assumed to be periodic or the initial data is compactly supported. In regions where the solution is smooth (i.e., some entropy functional of \cref{eq:hyp} satisfies a conservation equation), the system is reversible in the sense that it is invariant under the transformation $t \to -t$. However, in the vicinity of a shock, this time-symmetry property no longer holds as the entropy admissibility criterion induces a time irreversibility condition on the solution \citep{Dafermos2010_4}. This inconsistency can therefore be directly used to form a nonlinear stabilization technique for high-order numerical methods in the vicinity of discontinuous features. In particular, one can utilize a parabolic regularization of the conservation law, given in the form of 
\begin{equation}\label{eq:vp}
    \partial_t \mathbf{u} + \boldsymbol{\nabla}{\cdot} \mathbf{F}(\mathbf{u}) = \boldsymbol{\nabla}{\cdot} \left(\boldsymbol{\mu} \nabla \mathbf{u} \right),
\end{equation}
where $\boldsymbol{\mu}$ is some vector-valued artificial viscosity. The underlying physical mechanism for the proposed technique is closely related to the entropy viscosity approach of \citet{Guermond2011} where the amount of artificial viscosity is dependent on the amount of entropy production of the system. However, in contrast to solving an entropy residual equation to form a \emph{scalar-valued} viscosity, we propose an nonlinear artificial viscosity approach that can be readily implemented in a variety of numerical schemes without any additional governing equations, applies to systems without known entropy-flux pairs, and results in a \emph{vector-valued} viscosity. The component-wise decoupling of the viscous terms offers significant benefits for applications in nonlinear hyperbolic conservation laws, particularly for systems in which some components may exhibit discontinuities while others do not. Furthermore, although the proposed techniques are utilized to directly form an artificial viscosity term, the general approach is widely applicable to a variety of shock capturing methods in which an indicator metric is required.

The construction of this artificial viscosity is performed through an estimation of the reversibility of the system. Without loss of generality, we consider a forward temporal update with respect to a positive time step $\Delta t$ in terms of an explicit Euler approximation as 
\begin{equation}
    \overline{\mathbf{u}}^{n+1} = \mathbf{u}^{n} - \Delta t\left(\boldsymbol{\nabla}{\cdot} \mathbf{F}^{\, +} \left(\mathbf{u}^{n}\right)\right),
\end{equation}
where $\boldsymbol{\nabla}{\cdot} \mathbf{F}^{\, +}$ is some \emph{upwind}-biased approximation of the divergence of the flux from an arbitrary numerical scheme. The reversed state is then approximated from this state using a backward time step $-\Delta t$ as
\begin{equation}
    \mathbf{u}^{*} = \overline{\mathbf{u}}^{n+1} + \Delta t\left(\boldsymbol{\nabla}{\cdot} \mathbf{F}^{\, -}\left(\overline{\mathbf{u}}^{n+1}\right)\right),
\end{equation}
where $\boldsymbol{\nabla}{\cdot} \mathbf{F}^{\, -}$ is now a \emph{downwind}-biased approximation of the divergence of the flux (i.e., upwind with respect to the reverse time $-t$). Note that the time step and temporal integration scheme used to calculate the intermediate states does not necessarily have to be identical to the ones used for the overarching system which allows for a portion of the computation to be embedded within the overarching solver. 

The defect between the original state and the reversed state is then defined as
\begin{equation}
    \boldsymbol{\mathcal D} (\mathbf{u}^n) = |\mathbf{u}^{*} - \mathbf{u}^{n}|,
\end{equation}
which can be used as \emph{component-wise} indicator for the irreversibility of the system. In smooth (i.e., reversible) portions of the flow, this defect should tend to zero. Therefore, the defect would scale with the approximation error $\mathcal O (h^{p}\Delta t^{q+1})$, where $h$ is some characteristic local mesh scale and $p$ and $q$ are the spatial and temporal orders of the scheme, respectively. We posit, and show in \cref{sec:analysis}, that this term is expected to scale as $\mathcal O (h^{-1} \Delta t)$ in the vicinity of a shock. 

The artificial viscosity is then defined as 
\begin{equation}\label{eq:av}
    \boldsymbol{\mu}^n = c_\mu\frac{h^2}{\Delta t} \frac{\boldsymbol{\mathcal D} (\mathbf{u}^n)}{\Delta \mathbf{u}^n},
\end{equation}
where $c_\mu$ is a tunable constant and $\Delta \mathbf{u}^n$ is some normalization term generally taken as 
\begin{equation}
    \Delta \mathbf{u}^n =  \underset{\mathbf{x} \in \Omega}{\mathrm{max}} \left[\mathbf{u} (\mathbf{x}, \tau) \right] - \underset{\mathbf{x} \in \Omega}{\mathrm{min}} \left[\mathbf{u} (\mathbf{x}, \tau) \right]
\end{equation}
for $\tau \in [t^n, t^{n+1}]$. With these scaling factors, it can be seen that the components of $\boldsymbol{\mu}$ take the dimension of viscosity and this viscosity is of $\mathcal O (h)$ in the vicinity of a shock. We hereafter refer to this artificial viscosity as the time-reversibility (TR) viscosity in this work.

%% file: analysis.tex
\section{Analysis}
\label{sec:analysis}
Although physical solutions of the governing conservation laws should abide by the time-reversibility conditions, it is not immediately evident whether the numerical approximations of these conservation laws possess the same properties. To analyze the numerical mechanisms of time-reversibility, we considered the cases of scalar one-dimensional linear and nonlinear conservation laws of the form
\begin{equation}
    \partial_t u + \partial_x f(u) = 0.
\end{equation}
To simplify the analysis and show some useful properties in the limiting case, we consider a first-order spatial discretization on a uniform grid using an explicit Euler temporal integration scheme, given as
\begin{equation}
    u_i^{n+1} = u_i^n - \frac{\Delta t}{\Delta x} \left(f_{i+\shalf}^n - f_{i-\shalf}^n \right),
\end{equation}
where the subscript $i$ denotes the element index, the superscript $n$ denotes the temporal step (i.e., $t = n \Delta t$), $\Delta x$ denotes the mesh spacing (i.e., $x = i\Delta x$), and $f_{i-\shalf}^n/f_{i+\shalf}^n$ denote the left/right interface fluxes. This formulation is recovered for first-order finite volume and discontinuous Galerkin/flux reconstruction \citep{Hesthaven2008, Huynh2007} formulations on one-dimensional uniform grids. The reversed state is similarly computed as 
\begin{equation}
    u_i^{*} = u_i^{n+1} + \frac{\Delta t}{\Delta x} \left(f_{i+\shalf}^{n+1} - f_{i-\shalf}^{n+1} \right).
\end{equation}

The temporal integration scheme is assumed to abide by some standard CFL condition (i.e., $\lambda \frac{\Delta t}{\Delta x} < 1$ for some maximum wavespeed $\lambda$). We consider a discontinuous initial condition in the form of a step function, given as
\begin{equation}
\begin{cases}
     u_i^n = u_L, \quad \mathrm{for}\ i \leq 1\\
     u_i^n = u_R, \quad \mathrm{else},
\end{cases}
\end{equation}
where the magnitude of the discontinuity is given as
\begin{equation}
    \delta u = u_L - u_R > 0.
\end{equation}
The following analysis will primarily focus on the behavior of the TR viscosity in the vicinity of the discontinuity (i.e., $i=1,2$).

\subsection{Linear Transport}
The case of linear advection is first considered, where the flux is given as
\begin{equation}
    f(u) = \lambda u,
\end{equation}
for some constant advection velocity $\lambda > 0$. For a positive value of $\lambda$, the upwind-biased flux is simply the flux of the upwind state, i.e.,
\begin{equation}
    f_{i+\shalf}^n = f_i^n = \lambda u_i^n.
\end{equation}
From this, the intermediate states can be computed as
\begin{equation}
\begin{cases}
    u_0^{n+1} = u_0^{n} - \frac{\Delta t}{\Delta x} \left(f_{0}^n - f_{-1}^n \right) = u_0^{n}, \\
    
    u_1^{n+1} = u_1^n - \frac{\Delta t}{\Delta x} \left(f_{1}^n - f_{0}^n \right) = u_1^n - \lambda\frac{\Delta t}{\Delta x} \left(u_{1}^n - u_{0}^n \right)  = u_1^n,\\
    
    u_2^{n+1} = u_2^n - \frac{\Delta t}{\Delta x} \left(f_{2}^n - f_{1}^n \right) = u_2^n - \lambda\frac{\Delta t}{\Delta x} \left(u_{2}^n - u_{1}^n \right),\\
    
    u_3^{n+1} = u_3^n - \frac{\Delta t}{\Delta x} \left(f_{3}^n - f_{2}^n \right) = u_3^n - \lambda\frac{\Delta t}{\Delta x} \left(u_{3}^n - u_{2}^n \right)  = u_3^n.\\
    
\end{cases}
\end{equation}
To compute the reversed state, the downwind-biased flux is used, given as
\begin{equation}
    f_{i+\shalf}^{n+1} = f_{i+1}^{n+1} = \lambda u_{i+1}^{n+1},
\end{equation}
from which the reversed states can be computed as 
\begin{equation}
\begin{cases}
    u_1^{*} = u_1^{n+1} + \lambda\frac{\Delta t}{\Delta x} \left(u_{2}^{n+1} - u_{1}^{n+1} \right),\\
    
    u_2^{*} = u_2^{n+1} + \lambda\frac{\Delta t}{\Delta x} \left(u_{3}^{n+1} - u_{2}^{n+1} \right).
    
\end{cases}
\end{equation}

By substituting in the values of the intermediate states, the solution defect to the left and right of the discontinuity can then be given as 
\begin{align}
  \begin{aligned}
    \mathcal D (u_1) &=  \left |u_1^{n} - u_1^{*} \right |=  \left |u_1^n - \left( u_1^{n+1} + \lambda\frac{\Delta t}{\Delta x} \left(u_{2}^{n+1} - u_{1}^{n+1} \right) \right) \right |\\&=\left | \lambda\frac{\Delta t}{\Delta x} \left(u_2^n - \lambda\frac{\Delta t}{\Delta x} \left(u_{2}^n - u_{1}^n \right) - u_{1}^{n} \right)\right | = \left(1 - \lambda\frac{\Delta t}{\Delta x} \right)  \left( \lambda\frac{\Delta t}{\Delta x} \right)\left | \delta u\right | ,
\intertext{and}
    \mathcal D (u_2)& = \left | u_2^{n} - u_2^{*} \right | =  \left | u_2^{n} - \left( u_2^{n+1} + c \frac{\Delta t}{\Delta x} \left(u_{3}^{n+1} - u_{2}^{n+1} \right) \right)\right | \\&= \left |-\lambda\frac{\Delta t}{\Delta x} \left(u_{2}^n - u_{1}^n \right) + c \frac{\Delta t}{\Delta x} \left(u_3^n - \left(u_2^n - \lambda\frac{\Delta t}{\Delta x} \left(u_{2}^n - u_{1}^n \right) \right) \right)\right | = \left(1 - \lambda\frac{\Delta t}{\Delta x} \right)  \left( \lambda\frac{\Delta t}{\Delta x} \right)\left | \delta u\right | = \mathcal D (u_1),
  \end{aligned}
\end{align} 
respectively. Due to the symmetry of the forward/backward temporal integration, it can be seen that the defect is identical to the left and right of the discontinuity. From \cref{eq:av}, the TR viscosity can then be computed as 
\begin{equation}
    \mu = c_\mu \frac{\Delta x^2}{\Delta t} \frac{D(u_1)}{\left | \delta u \right |} = c_\mu ( \lambda  \Delta x - \lambda^2 \Delta t).
\end{equation}
In the limit as $\Delta t \to 0$, the TR viscosity in the vicinity of the discontinuity reduces to a particularly familiar form,
\begin{equation}
    \lim_{\Delta t \to 0} \mu = c_\mu \lambda  \Delta x.
\end{equation}
By setting $c_\mu = \shalf$, it can be seen that this corresponds to a `monotone' viscosity, where a second-order central scheme augmented with viscous dissipation exactly recovers a monotonic first-order upwind scheme. 

\subsection{Nonlinear Transport}
The analysis is similarly extended to the nonlinear case through the Burgers equation, where the flux is given as 
\begin{equation}
    f(u) = \frac{1}{2}u^2.
\end{equation}
For simplicity, we assume that $u_L, u_R > 0$, such that an upwind-biased flux can be simply computed as
\begin{equation}
    f_{i+\shalf}^n = f_i^n = \frac{1}{2}(u_i^n)^2.
\end{equation}
Similar conclusions can be drawn for the more general case utilizing a Rusanov-type upwind flux \citep{Rusanov1962}, although with significantly more algebraic complexity. As with the linear case, the forward states can be computed as 
\begin{equation}
\begin{cases}
    u_0^{n+1} = u_0^{n+1}, \\
    
    u_1^{n+1} = u_1^n - \frac{\Delta t}{2\Delta x} \left(f_{1}^n - f_{0}^n \right) = u_1^n - \frac{\Delta t}{2\Delta x} \left((u_1^n)^2- (u_0^n)^2 \right)  = u_1^n,\\
    
    u_2^{n+1} = u_2^n - \frac{\Delta t}{2\Delta x} \left(f_{2}^n - f_{1}^n \right) = u_2^n - \frac{\Delta t}{2\Delta x} \left((u_2^n)^2- (u_1^n)^2 \right),\\
    
    u_3^{n+1} = u_3^n - \frac{\Delta t}{2\Delta x} \left(f_{3}^n - f_{2}^n \right) = u_3^n - \frac{\Delta t}{2\Delta x} \left((u_3^n)^2- (u_2^n)^2 \right)  = u_3^n,\\
    
\end{cases}
\end{equation}
from which the reversed states can be computed as
\begin{equation}
\begin{cases}
    u_1^{*} = u_1^{n+1} + \frac{\Delta t}{2\Delta x} \left((u_{2}^{n+1})^2 - (u_{1}^{n+1})^2 \right),\\
    
    u_2^{*} = u_2^{n+1} +  \frac{\Delta t}{2\Delta x} \left((u_{3}^{n+1})^2 - (u_{2}^{n+1})^2 \right),\\
\end{cases}
\end{equation}
utilizing a downwind-biased flux,

\begin{equation}
    f_{i+\shalf}^{n+1} = f_{i+1}^{n+1} = \frac{1}{2}(u_{i+1}^{n+1})^2.
\end{equation}
Let $\zeta$ be defined as 
\begin{equation}
    \zeta = u_L^2 - u_R^2 > 0.
\end{equation}
The defect can then be expressed as 
\begin{align}
  \begin{aligned}
    \mathcal D (u_1) &= \left | u_1^{n} - u_1^{*}\right | = \left | u_1^n - \left(u_1^{n+1} + \frac{\Delta t}{2\Delta x} \left((u_{2}^{n+1})^2 - (u_{1}^{n+1})^2 \right) \right)\right | =\left |  -\frac{\Delta t}{2\Delta x} \left((u_{2}^{n+1})^2 - (u_{1}^{n+1})^2 \right)\right | \\&=\left |  -\frac{\Delta t}{2\Delta x} \left(\left(u_2^n - \frac{\Delta t}{2\Delta x} \left((u_2^n)^2- (u_1^n)^2 \right)\right)^2 - (u_{1}^{n})^2 \right) \right | = \left |  \frac{\Delta t}{2\Delta x} 
    \left( \zeta -  2 \left( \frac{\Delta t}{2\Delta x} \right) u_2^n \zeta + \left( \frac{\Delta t}{2\Delta x} \right)^2 \zeta ^2 \right) \right |,
  \end{aligned}
\end{align} 
and it can be shown that the defect remains symmetric in the nonlinear case,
\begin{equation}
    \mathcal D (u_2) = \mathcal D (u_1).
\end{equation}
With the addition of the scaling factors, the magnitude of the TR viscosity can be computed as
\begin{equation}
    \mu = c_\mu \frac{\Delta x^2}{\Delta t} \frac{D(u_1)}{\left | \delta u \right |} = c_\mu  \frac{\Delta x }{2 | \delta u |}\left| \zeta -  2 \left( \frac{\Delta t}{\Delta x} \right) u_2^n \zeta + \left( \frac{\Delta t}{\Delta x} \right)^2 \zeta ^2 \right|.
\end{equation}
When taken in the limit as $\Delta t \to 0$, the viscosity takes on a similar `monotone' form as with the linear case, 
\begin{equation}
    \lim_{\Delta t \to 0} \mu = c_\mu  \left | \frac{\zeta}{2\delta u}\right |\Delta x  = c_\mu \left | \frac{u_L^2 - u_R^2}{2(u_L - u_R)}\right | \Delta x = c_\mu \left | \frac{1}{2}(u_L + u_R) \right |\Delta x =c_\mu  \lambda\Delta x,
\end{equation}
where $\lambda$ is now the shock speed.

%% file: implementation.tex
\section{Implementation}\label{sec:implementation}
The proposed TR viscosity method was implemented within the framework of a high-order discontinuous spectral element method for unstructured meshes in multiple dimensions. In this section, we outline details regarding the implementation of the proposed techniques.

\subsection{Discretization}
The governing systems of equations were discretized using the nodal discontinuous Galerkin (DG) \citep{Hesthaven2008} method recovered via the flux reconstruction (FR) approach of \citet{Huynh2007}. In this approach, the domain $\Omega$ is partitioned into $N_e$ elements $\Omega_k$ such that $\Omega = \bigcup_{N_e}\Omega_k$ and $\Omega_i\cap\Omega_j=\emptyset$ for $i\neq j$. The approximate solution $\mathbf{u} (\mathbf{x})$ within each element $\Omega_k$ is represented via a nodal approximation of the form
\begin{equation}
    \mathbf{u} (\mathbf{x}) = \sum_{i = 0}^{N-1} \mathbf{u} (\mathbf{x}_i^u) {\phi}_i (\mathbf{x}),
\end{equation}
where $\mathbf{x}_i^u \ \forall \ i \in \{0,..., N-1\}$ is a set of $N$ solution nodes and ${\phi}_i (\mathbf{x})$ is a set of nodal basis functions with the property ${\phi}_i (\mathbf{x}_j^u) = \delta_{ij}$. The order of the approximation, represented by $\mathbb P_p$ for some order $p$, is defined as the maximal order of $\mathbf{u} (\mathbf{x})$. Communication between elements is performed via the element interfaces $\partial \Omega$. We use the notation $\mathbf{p}^{\partial \Omega}$ to denote the interpolation of the function $\mathbf{p}(\mathbf{x})$ onto the set of $M$ interface points $\mathbf{x}_i^I \in \partial \Omega \ \forall \ i \in \{0,..., M-1\}$. Given an interface point $\mathbf{x}_i^I$, let $\mathbf{p}^{\partial \Omega^+}_i$ denote the value of $\mathbf{p}(\mathbf{x})$ evaluated from the element of interest and let $\mathbf{p}^{\partial \Omega^-}_i$ denote the value of $\mathbf{p}(\mathbf{x})$ evaluated from the interface-adjacent element. 

The formation of the approximate flux is split into two components, a nonlinear inviscid flux $\mathbf{F} \left (\mathbf{u} \right)$ corresponding to the flux of the \cref{eq:hyp} and a linear viscous flux $\mathbf{G} \left (\mathbf{u}, \mathbf{v} \right)$ corresponding to the parabolic regularization term in \cref{eq:vp}. An auxiliary variable $\mathbf{v}$ is introduced to represent the corresponding gradients of the solution such that $\mathbf{v} \approx \nabla \mathbf{u}$. Within the flux reconstruction methodology, the approximate flux is formed through an element-local approximation which is followed by a correction to account for interface contributions. 

Beginning with the inviscid flux, a discontinuous flux $\hat{\mathbf{f}} (\mathbf{x})$ is computed through a collocation projection of the flux function onto the solution nodes.
\begin{equation}
    \hat{\mathbf{f}} (\mathbf{x}) = \sum_{i = 0}^{N-1} \mathbf{F} \left (\mathbf{u} (\mathbf{x}_i^u) \right ) {\phi}_i (\mathbf{x}).
\end{equation}
This discontinuous flux is then corrected with the flux contributions at the element interfaces $\partial \Omega$. The approximate solution is interpolated to the interfaces to yield solution pairs for each interface flux point $\mathbf{x}_i^I$, denoted as $\{\mathbf{u}^{\partial \Omega^-}_i, \mathbf{u}^{\partial \Omega^+}_i\}$. These pairs, along with the corresponding outward facing normals $\mathbf{n}_i$, are used to form common interface inviscid fluxes ${\mathbf{f}}^{I}_i$ via an approximate or exact Riemann solver using approaches such as that of \citet{Rusanov1962} and \citet{Roe1981}. The corrected flux $\mathbf{f} (\mathbf{x})$ is then defined as
\begin{equation}
    \mathbf{f} (\mathbf{x}) = \hat{\mathbf{f}} (\mathbf{x}) +  \sum_{i = 0}^{M-1}\left(\mathbf{f}^{I}_i - \hat{\mathbf{f}}^{\partial \Omega^+}_i \right) \mathbf{h}_i(\mathbf{x})
\end{equation}
using a set of correction functions $\mathbf{h}_i(\mathbf{x})$. The correction functions \citep{Castonguay2011, Trojak2021} have the property that
\begin{equation*}
    \mathbf{n}_i \cdot \mathbf{h}_j(\mathbf{x}_i) = \delta_{ij} \quad \quad \mathrm{and} \quad \quad \sum_{i=0}^{M-1} \mathbf{h}_i(\mathbf{x}) \in \mathrm{RT}_p,
\end{equation*}
where $\mathrm{RT}_p$ is the Raviart--Thomas space of order $p$ \citep{Raviart1977}. In this work, the correction functions are chosen to recover the nodal discontinuous Galerkin scheme, as presented by \citet{Huynh2007}.

A similar method is performed with the viscous flux, although with a second set of corrections for the gradient terms. A discontinuous solution gradient is formed as $\hat{\mathbf{v}} = \nabla \mathbf{u}$, from which the corrected solution gradient can be defined as
\begin{equation}
    \mathbf{v} (\mathbf{x}) = \hat{\mathbf{v}} (\mathbf{x}) +  \sum_{i = 0}^{M-1}\left(\mathbf{u}^{I}_i - \mathbf{u}^{\partial \Omega^+}_i \right) \nabla \mathbf{h}_i(\mathbf{x}).
\end{equation}
The common interface solution $\mathbf{u}^{I}_i$ can be calculated as a centered mean or using approaches such as that of \citet{Bassi2000}. A discontinuous viscous flux can then be formed as 
\begin{equation}
    \hat{\mathbf{g}} (\mathbf{x}) = \sum_{i = 0}^{N-1} \mathbf{G} \left (\mathbf{u} (\mathbf{x}_i^u), \mathbf{v} (\mathbf{x}_i^u) \right ) {\phi}_i (\mathbf{x}),
\end{equation}
from which the corrected viscous flux can be formed as 
\begin{equation}
    \mathbf{g} (\mathbf{x}) = \hat{\mathbf{g}} (\mathbf{x}) +  \sum_{i = 0}^{M-1}\left(\mathbf{g}^{I}_i - \hat{\mathbf{g}}^{\partial \Omega^+}_i \right) \mathbf{h}_i(\mathbf{x}).
\end{equation}
The common interface viscous flux $\mathbf{g}^{I}_i$ can be similarly computed as with the common interface solution.
From this, the purely hyperbolic and the viscous regularized semidiscretizations can then be formed as
\begin{align}
    \partial_t \mathbf{u} &= -\boldsymbol{\nabla}{\cdot}\mathbf{f} (\mathbf{x}),\\
    \partial_t \mathbf{u} &= -\boldsymbol{\nabla}{\cdot}\mathbf{f} (\mathbf{x}) + \boldsymbol{\nabla}{\cdot}\mathbf{g} (\mathbf{x}),
\end{align}
respectively.

\subsection{Numerical Framework}
The proposed approach was implemented within the PyFR software package \citep{Witherden2014}, an unstructured flux reconstruction solver that can target multiple compute architectures including CPUs and GPUs. The governing equations were solved with the FR methodology using a Rusanov-type \citep{Rusanov1962} Riemann solver with Davis wavespeeds \citep{Davis1988} for the inviscid fluxes, the BR1 method of \citet{Bassi1997} for the viscous fluxes, and an explicit fourth-order Runge--Kutta (RK) scheme for temporal integration. For the calculation of the TR viscosity, a forward Euler approximation was used at each RK stage, and downwinding was performed by negating the sign of the diffusive component in the Riemann solver. The local mesh scale $h_k$ for each element $\Omega_k$ was calculated as
\begin{equation}
    h_k^d = \int_{\Omega_k} \mathrm{d} \mathbf{V},
\end{equation}
for a given spatial dimension $d$. Furthermore, in this work, an element-wise constant TR viscosity was used, taken as the integrated mean value over the element. Potential improvements in performance are possible by increasing the regularity of the viscosity profile, and the reader is referred to the works of \citet{Barter2007} and \citet{Glaubitz2018} for details. For the following problems, the free parameter $c_{\mu}$ was set to a scalar value, and some tuning of this parameter was performed on a problem-dependent basis except in the cases where the resolution/order was modified for the same problem. The acceptable range for $c_\mu$ was found to be on the order of 1-10.

Due to the use of an explicit time integration scheme, transients in the TR viscosity could severely limit the maximum acceptable time step. To remedy this, we set a maximum allowable global viscosity value, 
\begin{equation}
    \mu_{\mathrm{max}} = c_{\mathrm{max}} \hat{\lambda} h_{\mathrm{max}} ,
\end{equation}
where $c_{\mathrm{max}} \approx 100$ is some arbitrarily large constant, $ \hat{\lambda}$ is an \textit{a priori} estimate of the maximum wavespeed in the domain across the time range of the simulation, and $h_{\mathrm{max}}$ is the largest characteristic mesh size in the domain. The maximum viscosity value is computed once prior to the simulation and held constant throughout. This arbitrarily large bound is used purposefully as opposed to a more appropriate local viscosity bound as this ensures that the relation between the artificial viscosity magnitude and the time-reversibility metric remains proportional and that the effects of the viscosity limiting are purely to mitigate time step restrictions. An alternative approach is to limit the viscosity by a more appropriate local bound (see \citet{Guermond2011}), but this would introduce another tunable parameter and would obfuscate the ability of the proposed approach to predict the appropriate amount of required artificial viscosity. For the one-dimensional test cases, the viscous terms were treated implicitly, and no differences in the results were observed with and without this limiting. 

%% file: results.tex
\section{Results}\label{sec:results}
The accuracy and robustness of the TR viscosity method was evaluated on a variety of nonlinear hyperbolic conservation laws that exhibit discontinuous features. To show the utility of the approach, the numerical experiments were generally focused on very high order discretizations.  
\subsection{Nonlinear Transport}

\subsubsection{Burgers Equation}
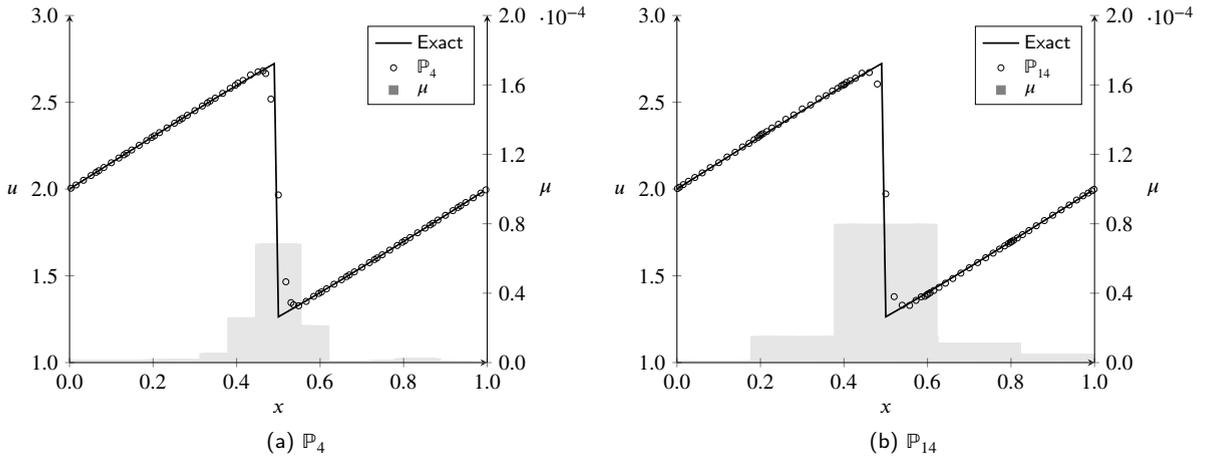
\begin{figure}[tbhp]
        \centering        
        \subfloat[$\mathbb P_4$]{
         \adjustbox{width=0.48\linewidth,valign=b}{\input{./figs/burgers_15}}}
        \subfloat[$\mathbb P_{14}$]{
         \adjustbox{width=0.48\linewidth,valign=b}{\input{./figs/burgers_5}}}
        \caption{\label{fig:burgers} Solution of the Burgers equation at $t = 0.5$ computed using a $\mathbb P_{4}$ FR scheme with $15$ elements (left) and $\mathbb P_{14}$ FR scheme with $5$ elements (right). Circles denote the individual solution nodes.}
    \end{figure}

The preliminary evaluation of the proposed method was performed on the Burgers equation, written as 
    \begin{equation}
        \partial_t u + \partial_x \left (\half u^2 \right) = 0.
    \end{equation}
An initially smooth profile was imposed on the periodic domain $\Omega = [0,1]$, given by the initial condition
    \begin{equation*}
        u(x, 0) = \sin(2\pi x) + 2.
    \end{equation*}
To demonstrate the efficacy of the TR viscosity approach in the context of high-order methods, a very coarse mesh consisting of 5 elements was used and discretized with a $\mathbb{P}_{14}$ FR scheme and compared to a finer mesh consisting of 15 elements using a $\mathbb{P}_{4}$ FR scheme. The predicted solution at $t = 0.5$ is shown in \cref{fig:burgers} along with the distribution of the TR viscosity ($c_\mu = 1$). Even with the very high approximation degree, the shock within the middle element was resolved over a span of 4-5 solution points for the $\mathbb P_{14}$ scheme, a small portion of the total element size. Furthermore, no spurious oscillations were observed in the solution profile, and the TR viscosity distribution was focused in the vicinity of the discontinuity. In comparison, the lower-order $\mathbb P_{4}$ scheme was able to resolve the discontinuity over a similar distance, indicating that the proposed approach does not suffer with increasing approximation order. 

\subsubsection{Kurganov--Petrova--Popov Problem}
\begin{figure}
        \centering
        \subfloat[Solution]{\label{fig:kpp_u} \adjustbox{width=0.4\linewidth,valign=b}{\includegraphics[width=\textwidth]{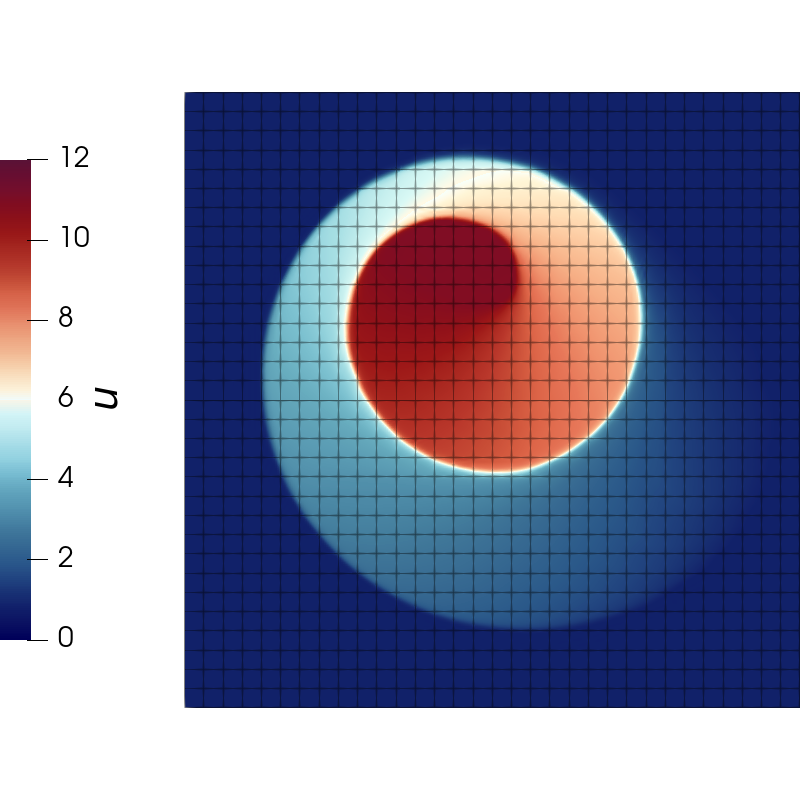}}}
        ~
        \subfloat[Viscosity]{\label{fig:kpp_v} \adjustbox{width=0.4\linewidth,valign=b}{\includegraphics[width=\textwidth]{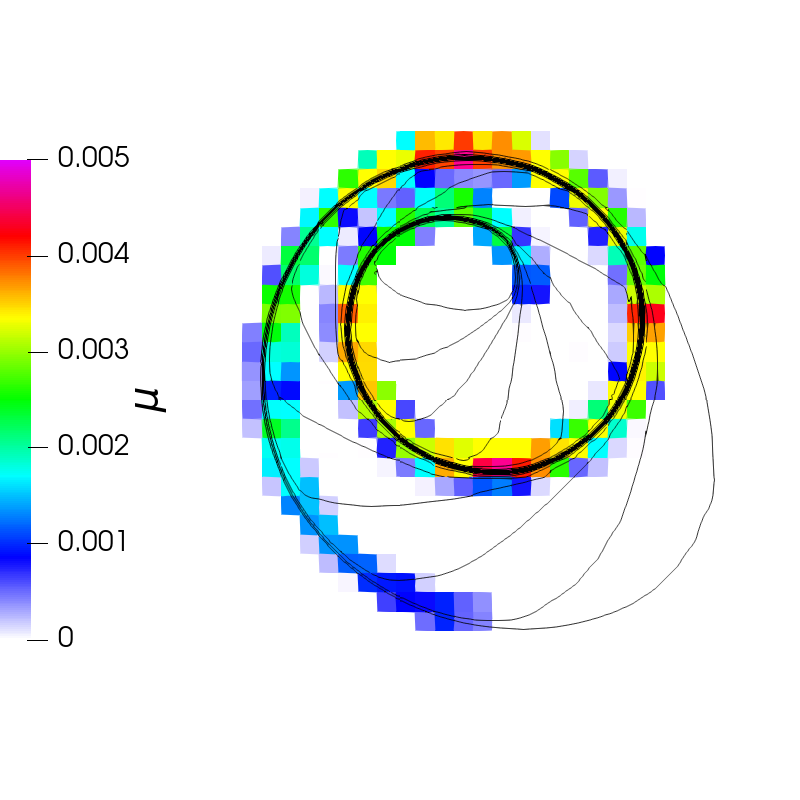}}}
        \newline
        \caption{\label{fig:kpp} Solution contours (left) and TR viscosity (right) for the KPP problem at $t = 1$ computed using a $\mathbb P_9$ FR scheme with $32^2$ elements. Solution contours are overlayed with the mesh and viscosity contours are overlayed with 20 equispaced isocontours of the solution. }
    \end{figure}

\begin{figure}
        \centering
        \subfloat[$N = 16^2$]{\adjustbox{width=0.24\linewidth,valign=b}{\includegraphics[width=\textwidth]{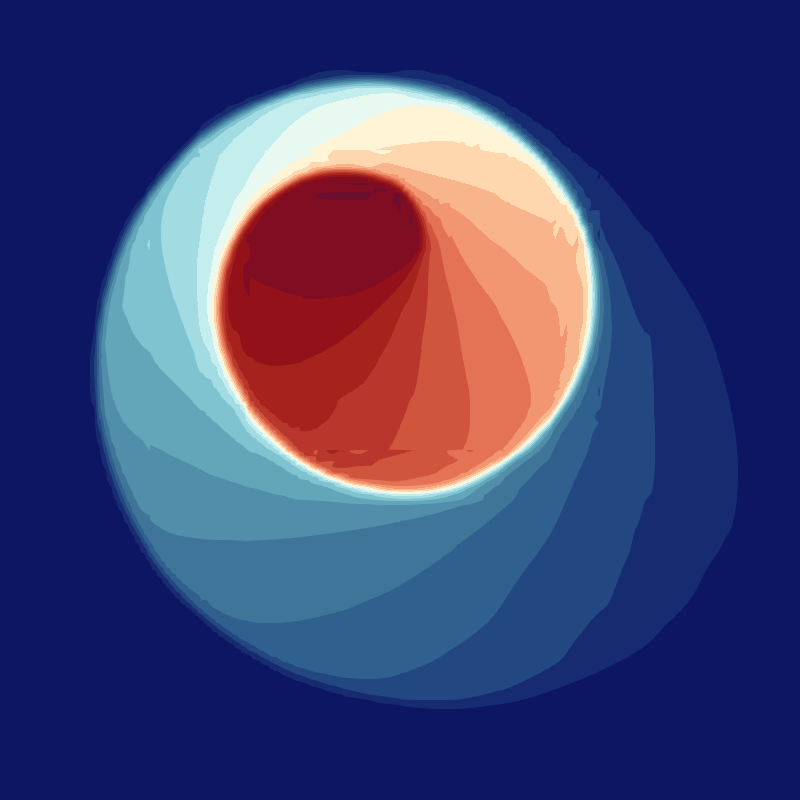}}}
        \hfill
        \subfloat[$N = 32^2$]{\adjustbox{width=0.24\linewidth,valign=b}{\includegraphics[width=\textwidth]{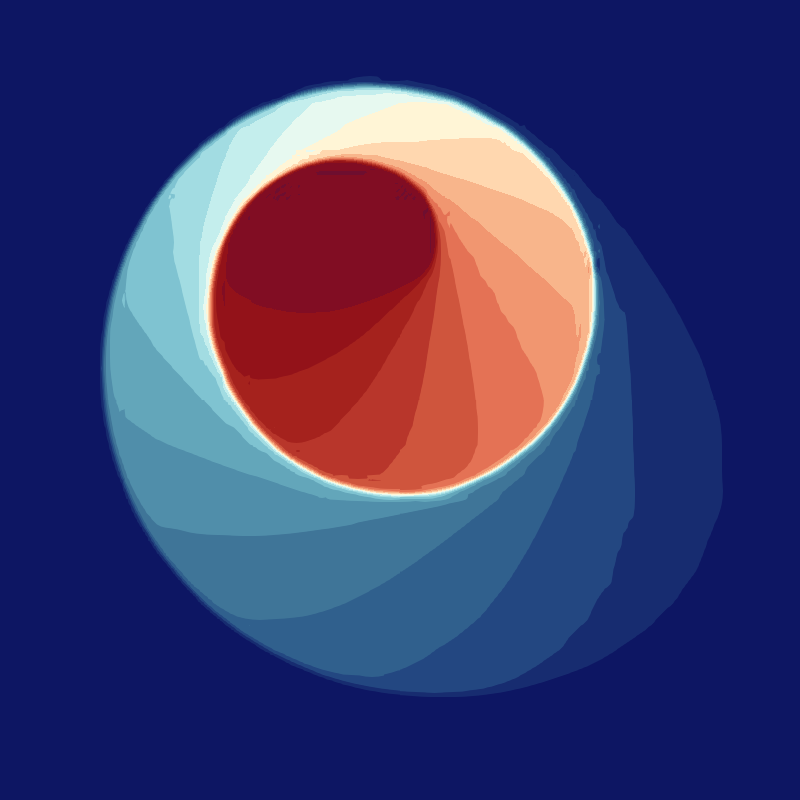}}}
        \hfill
        \subfloat[$N = 48^2$]{\adjustbox{width=0.24\linewidth,valign=b}{\includegraphics[width=\textwidth]{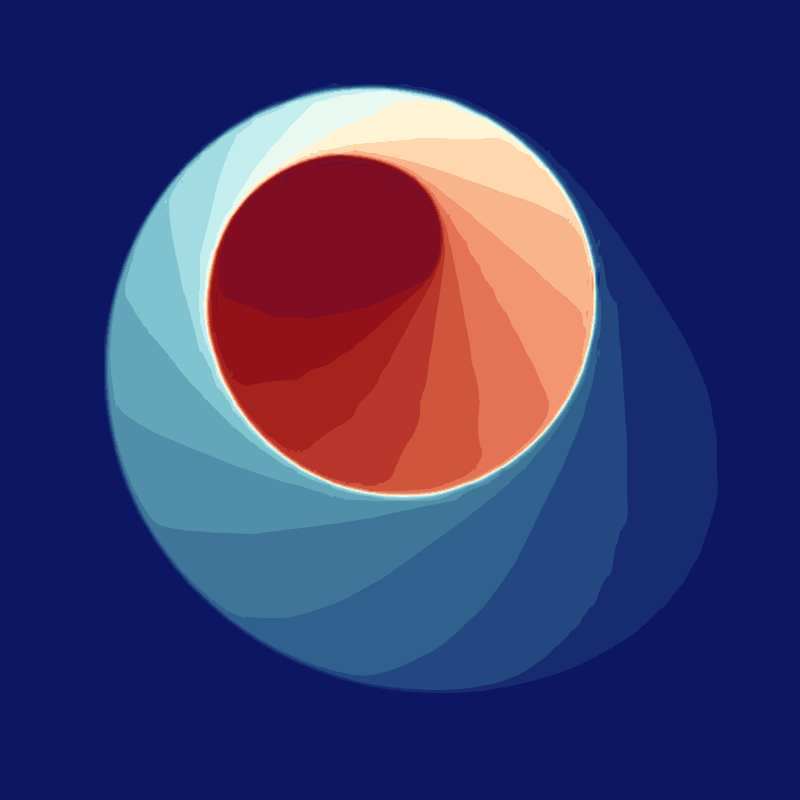}}}
        \hfill
        \subfloat[$N = 64^2$]{\adjustbox{width=0.24\linewidth,valign=b}{\includegraphics[width=\textwidth]{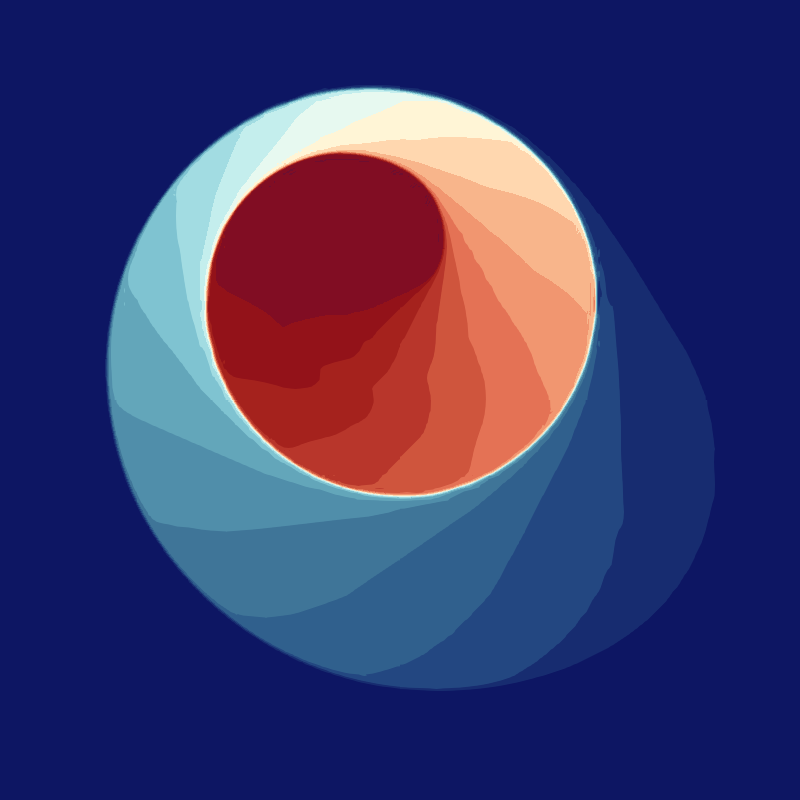}}}
        ~
        \newline
        \caption{\label{fig:kpp_conv} Solution contours for the KPP problem at $t = 1$ computed using a $\mathbb P_9$ FR scheme with $16^2$, $32^2$, $48^2$, and $64^2$ elements. }
    \end{figure}

    The extension to multi-dimensional scalar conservation laws was conducted through the Kurganov--Petrova--Popov (KPP) rotating wave problem \citep{Kurganov2007}. The governing conservation law is defined as
    \begin{equation}
        \partial_t u + \boldsymbol{\nabla}{\cdot} \left[ \sin u, \cos u \right]^T = 0,
    \end{equation}
    with the initial conditions defined by 
    \begin{equation*}
        u(\mathbf{x}, 0) = \begin{cases}
            3.5 \pi, &\mbox{if } x^2 + y^2 \leq 1, \\
            \frac{1}{4}\pi, &\mbox{else},
        \end{cases} 
    \end{equation*}
    on the periodic domain $\Omega = [-2, 2]^2$. The piecewise-constant initial data develops into a two-dimensional composite wave structure, and as a result, this problem poses a challenge to many high-resolution schemes with many schemes failing to converge to the proper entropy solution. The results of the TR viscosity method ($c_\mu = 6$) at $t=1$ using a $\mathbb P_9$ FR scheme with $32^2$ elements are shown in \cref{fig:kpp}. Contours of the solution show the expected wave structure of the problem, and the discontinuities in the wave are properly resolved with sub-element resolution without spurious oscillations even at a very high order. The distribution of the TR viscosity shows that its effects are localized around the discontinuities and roughly proportional to the magnitude of the discontinuity. To qualitatively verify that the predicted solution is converging with respect to mesh size, the solution contours as computed on an array of meshes with increasing resolution are shown in \cref{fig:kpp_conv}. For a given approximation order and fixed value of $c_\mu$, the results showed qualitative convergence to the expected wave structure of the problem and sharper resolution of the discontinuities with decreasing mesh size. 

\subsection{Euler Equations}
    The TR viscosity method was then assessed for the compressible Euler equations, written as 
    \begin{equation}
        \mathbf{u} = \begin{bmatrix}
                \rho \\ \boldsymbol{\rho v} \\ E
            \end{bmatrix}, \quad  \mathbf{F} = \begin{bmatrix}
                \boldsymbol{\rho v}\\
                \boldsymbol{\rho v}\otimes\mathbf{v} + P\mathbf{I}\\
            (E+P)\mathbf{v}
        \end{bmatrix},
    \end{equation}
    where $\rho$ is the density, $\boldsymbol{\rho v}$ is the momentum, $E$ is the total energy, $P = (\gamma-1)\left(E - \shalf\rho\mathbf{v}{\cdot}\mathbf{v}\right)$ is the pressure, and $\gamma = 1.4$ is the ratio of specific heat capacities. The symbol $\mathbf{I}$ denotes the identity matrix in $\mathbb{R}^{d\times d}$ and $\mathbf{v} = \boldsymbol{\rho v}/\rho$ denotes the velocity. For brevity, the solution is expressed in terms of a vector of primitive variables as $\mathbf{q}=[\rho,\mathbf{v},P]^T$.
\subsubsection{Sod Shock Tube}
    The Sod shock tube problem \citep{Sod1978} was used as an evaluation of the method in predicting the three main features of the Riemann problem: rarefaction waves, contact discontinuities, and shock waves. The problem is solved on the domain $\Omega=[0,1]$ with the initial conditions
    \begin{equation*}
        \mathbf{q}(x,0) = \begin{cases}
            \mathbf{q}_l, &\mbox{if } x\leqslant 0.5,\\
            \mathbf{q}_r, &\mbox{else},
        \end{cases} \quad \mathrm{given} \quad \mathbf{q}_l = \begin{bmatrix}
            1 \\ 0 \\ 1
        \end{bmatrix}, \quad \mathbf{q}_r = \begin{bmatrix}
            0.125 \\ 0 \\ 0.1
        \end{bmatrix}.
    \end{equation*}
    
\begin{figure}[tbhp]
        \centering        
        \subfloat[$\mathbb P_3$]{\label{fig:sod_d} \adjustbox{width=0.48\linewidth,valign=b}{\input{./figs/sod_density_p3}}}
        ~
        \subfloat[$\mathbb P_9$]{\label{fig:sod_v} \adjustbox{width=0.48\linewidth,valign=b}{\input{./figs/sod_density_p9}}}
        \newline
        \caption{\label{fig:sod} Density profile of the Sod shock tube at $t=0.2$ computed using a $\mathbb P_3$ (left) and $\mathbb P_9$ (right) FR scheme with 200 degrees of freedom.}
    \end{figure}
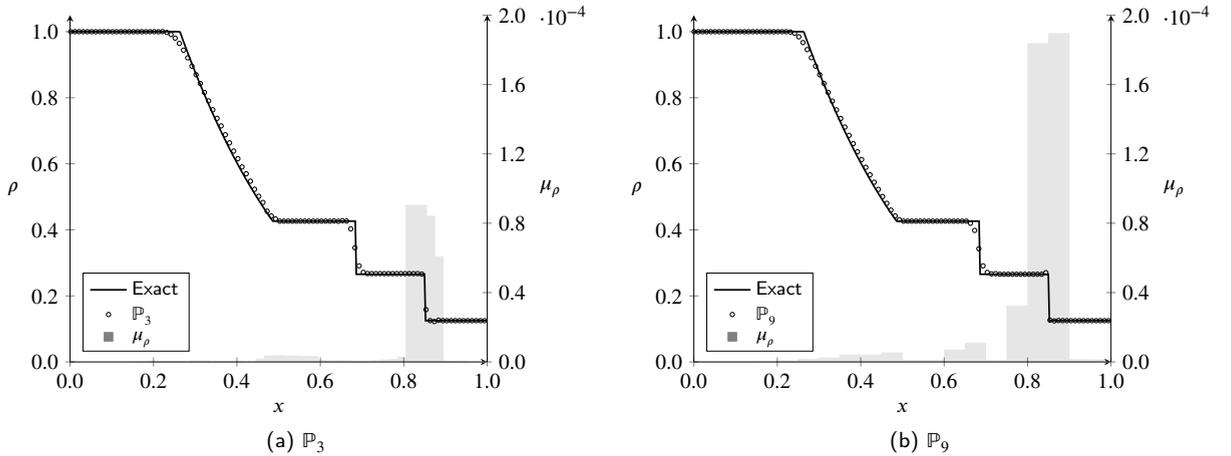

    The results of the shock tube at $t = 0.2$ using the TR viscosity method ($c_\mu = 5$) are presented in \cref{fig:sod}. The degrees of freedom were fixed at 200 for both the $\mathbb P_3$ and $\mathbb P_9$ schemes, corresponding to 50 and 20 elements, respectively, across the domain. For the relatively coarse resolution, the results showed excellent agreement with the exact solution, with good resolution of the shock front and contact discontinuity and minimal spurious oscillations. This resolution was maintained even when the order of the approximation was increased, such that the discontinuities were resolved well within the element. Furthermore, the distribution of the density component of the TR viscosity was localized around the shock front, and minimal additional dissipation was introduced around the rarefaction wave and contact discontinuity. 
    
    To quantitatively verify the convergence of the entropy filtering approach for discontinuous solutions, the convergence rates of the error against the exact solution were evaluated, specifically the point-mean $L^1$ norm of the density error, defined as
    \begin{equation}
        \epsilon_{\rho_1} = \frac{1}{M}\sum_{i=0}^{M-1} \left|\rho (x_i) - \rho_{\mathrm{exact}} (x_i) \right|,
    \end{equation}
    where $M$ is the total degrees of freedom. The convergence rates of the density error are shown in \cref{tab:sod_error} with respect to the number of elements $N$ for varying approximation orders. The expected first-order convergence rate was recovered for all approximation orders, and for a given number of degrees of freedom, a higher approximation order generally resulted in marginally lower error.

    \begin{figure}[tbhp]
        \centering
        \begin{tabular}{r | c c c c c c}\toprule
	        $N$ & $\mathbb{P}_2$ &$\mathbb{P}_3$ &$\mathbb{P}_4$ &$\mathbb{P}_5$ &$\mathbb{P}_6$ &$\mathbb{P}_7$ \\ \midrule
	        
            $20$ & \num{4.15e-03} & \num{2.79e-03} & \num{1.80e-03} & \num{1.52e-03} & \num{1.34e-03} & \num{1.01e-03} \\
            $40$ & \num{1.93e-03} & \num{1.29e-03} & \num{8.68e-04} & \num{7.75e-04} & \num{5.48e-04} & \num{5.87e-04} \\
            $80$ & \num{9.63e-04} & \num{5.86e-04} & \num{4.50e-04} & \num{4.14e-04} & \num{3.08e-04} & \num{3.13e-04} \\
            $160$ & \num{6.24e-04} & \num{3.43e-04} & \num{3.03e-04} & \num{2.10e-04} & \num{1.66e-04} & \num{1.77e-04} \\
            \midrule
        \textbf{RoC} & $0.92$& ${1.02}$& ${0.87}$& ${0.95}$& ${0.99}$& ${0.85}$\\ 
        \end{tabular}
        \captionof{table}{\label{tab:sod_error} $L_1$ norm of the density error for the Sod shock tube problem at varying orders. Rate of convergence shown on bottom.}
    \end{figure}

\subsubsection{Shu--Osher Problem}
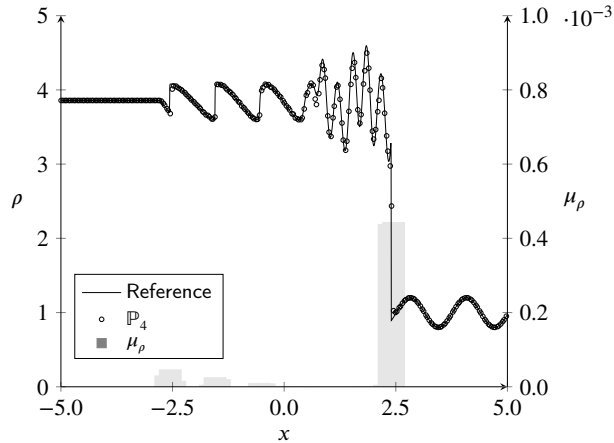
\begin{figure}[tbhp]
        \centering
        \adjustbox{width=0.5\linewidth,valign=b}{\input{./figs/shu_density}}
        \caption{\label{fig:shu} Density profile of the Shu-Osher problem at $t = 1.8$ computed using a $\mathbb P_4$ FR scheme with 100 elements.}
    \end{figure}

    The more complex problem of \citet{Shu1988} poses a challenge for schemes in their ability to accurately identify and resolve discontinuities in the presence of physical oscillations. This problem, solved on the domain $\Omega = [-5, 5]$ with the initial conditions
    \begin{equation*}
         \mathbf{q}(x,0) =  \begin{cases}
            \mathbf{q}_l, &\mbox{if } x\leqslant -4, \\
            \mathbf{q}_r, &\mbox{else},
        \end{cases} \quad \mathrm{given} \quad
        \mathbf{q}_l = \begin{bmatrix}
            3.857143 \\ 2.629369 \\ 10.333333
        \end{bmatrix}, \quad
        \mathbf{q}_r = \begin{bmatrix}
            1 + 0.2\sin{5x} \\ 0\\ 1
        \end{bmatrix},
    \end{equation*}
    consists of a shock front propagating through a sinusoidally-perturbed density field. This interaction yields instabilities in the flow field which may be erroneously identified as shocks and, as a result, excessively dissipated by certain schemes. The density profile as computed by the TR viscosity method ($c_\mu = 5$) using a $\mathbb P_4$ FR scheme with 100 elements is shown in \cref{fig:shu}. A reference solution was obtained via a highly-resolved exact Godunov-type solver \cite{Toro1997_4}. The results show good agreement with the reference profile, notably with the resolution of the shock front and trailing shocks as well as the instabilities in the density field. The viscosity distribution shows that the TR indicator identified the shock waves distinctly from the oscillatory density field and the magnitude of the viscosity was roughly proportional to the shock strength. 

\subsubsection{Isentropic Euler Vortex}

To verify that the TR viscosity approach does not deteriorate the accuracy for smooth solutions, the isentropic Euler vortex \citep{Shu1998} is used as the problem has an analytic solution that can be used to evaluate the convergence of the error. The initial conditions are given as 
    \begin{equation}
        \mathbf{q}(\mathbf{x}, 0)  = \begin{bmatrix}
            p^\frac{1}{\gamma} \\
            V_x + \frac{S}{2 \pi R} (y-y_0)\phi(r) \\[4pt]
            V_y - \frac{S}{2 \pi R} (x-x_0)\phi(r) \\[4pt]
            \frac{1}{\gamma M^2} \left(1 - \frac{S^2 M^2 (\gamma-1)} {8 \pi^2}\phi(r)^2\right)^\frac{\gamma}{\gamma-1}
        \end{bmatrix}, \quad \mathrm{where} \quad r = \|\mathbf{x}-\mathbf{x}_0\|_2 \quad \mathrm{and} \quad \phi(r) = \exp{\left(\frac{1-r^2}{2R^2}\right)}.
    \end{equation}
The characteristics of the vortex are defined by the parameters $S = 13.5$ denoting the strength of the vortex, $R = 1.5$ the radius, $V_x = 0$, $V_y = 1$ the advection velocities, and $M = 0.4$ the freestream Mach number. A uniform quadrilateral mesh on the domain $\Omega = [-10, 10]^2$ was generated and periodic boundary conditions were applied. The value of the $c_\mu$ parameter was set to 6 as this was the largest value used for the remaining two-dimensional test cases. 

After one convective time in which the vortex has returned to its original position, the $L^2$ norm of the density error was calculated for a series of grids of $N \times N$ elements. The sequences of the error for varying orders are shown in \cref{tab:icv_error}. The recovered rates of convergence met or exceeded the expected theoretical convergence rates of order $p{+}1$. The distributions of the TR viscosity components were examined at the final states to verify that their magnitudes were negligible. For the $\mathbb{P}_7$ approaches, the maximum values of the four components were [\num{2.6E-08}, \num{8.6E-08}, \num{7.4E-08}, \num{4.7E-07}] for $N = 20$ and [\num{8.1E-11}, \num{1.7E-10}, \num{1.4E-10}, \num{1.1E-09}] for $N = 40$. These observations are consistent with assumption that the TR viscosity should tend towards the approximation error for smooth solutions and should converge to (machine) zero in the limit of infinite resolution. 

    \begin{figure}[tbhp]
        \centering
        \begin{tabular}{r | c c c c c c}\toprule
	        $N$ & $\mathbb{P}_2$ &$\mathbb{P}_3$ &$\mathbb{P}_4$ &$\mathbb{P}_5$ &$\mathbb{P}_6$ &$\mathbb{P}_7$ \\ \midrule
	        
	        $20$ & - & - & - & - & \num{3.67E-05} & \num{6.34E-06} \\ 
	        
	        $25$ & - & - & \num{1.01E-03} & \num{1.39E-04} & \num{6.70E-06} & \num{1.44E-06}  \\
	        
            $33$ & \num{5.83E-02}  & \num{3.46E-03} & \num{2.06E-04} & \num{1.65E-05} & \num{6.70E-07}  & \num{9.48E-08}\\
            
            $40$ & \num{3.18E-02} & \num{1.33E-03} & \num{7.24E-05} & \num{3.98E-06} & \num{1.75E-07} & \num{1.64E-08} \\ 
            
            $50$ & \num{1.54E-02} & \num{4.40E-04} & \num{2.19E-05} & \num{7.99E-07} & - & -  \\
            
            $67$ & \num{5.66E-03} & \num{1.05E-04} & - & - & - & - \\\midrule
        \textbf{RoC} & $3.30$& ${4.93}$& ${5.53}$& ${7.44}$& ${7.79}$& ${8.74}$\\ 
        \end{tabular}
        \captionof{table}{\label{tab:icv_error} $L_2$ norm of the density error for the isentropic Euler vortex problem at varying orders. Rate of convergence shown on bottom.}
    \end{figure}

\subsubsection{Partially-Confined Explosion}
\begin{figure}[h!]
        \centering
        \adjustbox{width=0.4\linewidth,valign=b}{\input{./figs/exp_geo}}
        \caption{\label{fig:expgeo} Schematic of the partially-confined explosion geometry on the domain $[0,1]^2$. Solid lines denote slip wall boundaries and red line denotes the left/right state separation diaphragm.}
    \end{figure}
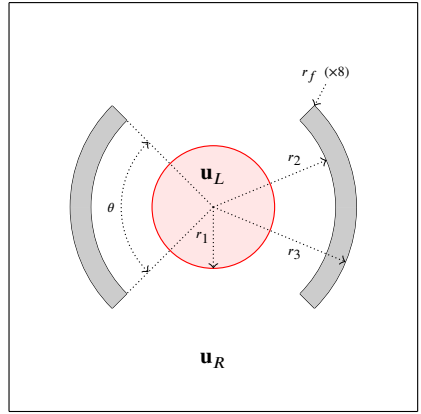

For an evaluation of the proposed method for more complex transient flows on two-dimensional unstructured grids, a novel problem setup, the partially-confined explosion, is introduced in this work. The problem, shown in \cref{fig:expgeo}, is defined on the domain $\Omega = [-1, 1]$ and consists of a compressed gas separated by a diaphragm placed within a circular region of radius $r_1 = 0.3$ centered at the origin. A circular confinement ring is placed in the region between $r_2 = 0.6$ and $r_3 = 0.7$ with a coverage angle of $\theta = 90^{\circ}$. The initial conditions are defined by 
    \begin{equation*}
        \mathbf{q}(x,0) = \begin{cases}
            \mathbf{q}_l, &\mbox{if } \| \mathbf{x} \|_2 \leqslant 0.3,\\
            \mathbf{q}_r, &\mbox{else},
        \end{cases} \quad \mathrm{given} \quad \mathbf{q}_l = \begin{bmatrix}
            1 \\ 0 \\ 1
        \end{bmatrix}, \quad \mathbf{q}_r = \begin{bmatrix}
            0.125 \\ 0 \\ 0.1
        \end{bmatrix},
    \end{equation*}
which yield an initial setup similar to a radially-symmetric Sod shock tube. The problem setup emulates a simplification of a lens-focused explosion, albeit with circular lens shape instead of a paraboloid lens shape. Two unstructured second-order triangular meshes of varying resolution were generated: a coarse mesh of $2\cdot10^5$ elements and a fine mesh of $8\cdot10^5$, corresponding to an average element edge length $\bar{h}$ of approximately 0.008 and 0.004, respectively. For all boundaries, slip-wall boundary condition were imposed. To alleviate numerical issues with regards to sharp corners, a fillet of radius $r_f = 0.01$ was applied to all corners of the confinement ring. 

    \begin{figure}[h!]
        \centering
        \subfloat[$t = 1$]{\label{fig:ant_0} \adjustbox{width=0.3\linewidth,valign=b}{\includegraphics[width=\textwidth]{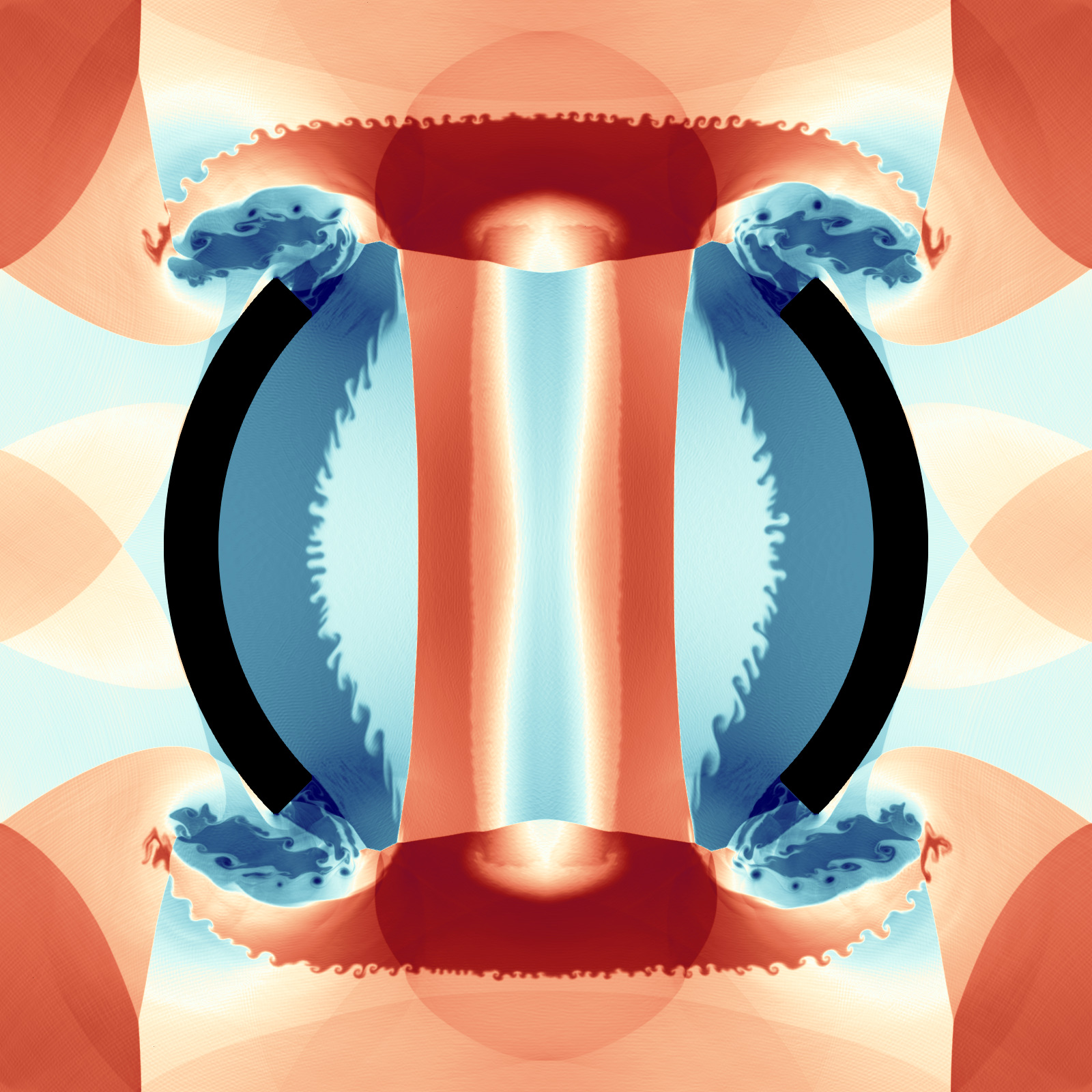}}}
        ~
        \subfloat[$t = 2$]{\label{fig:ant_1} \adjustbox{width=0.3\linewidth,valign=b}{\includegraphics[width=\textwidth]{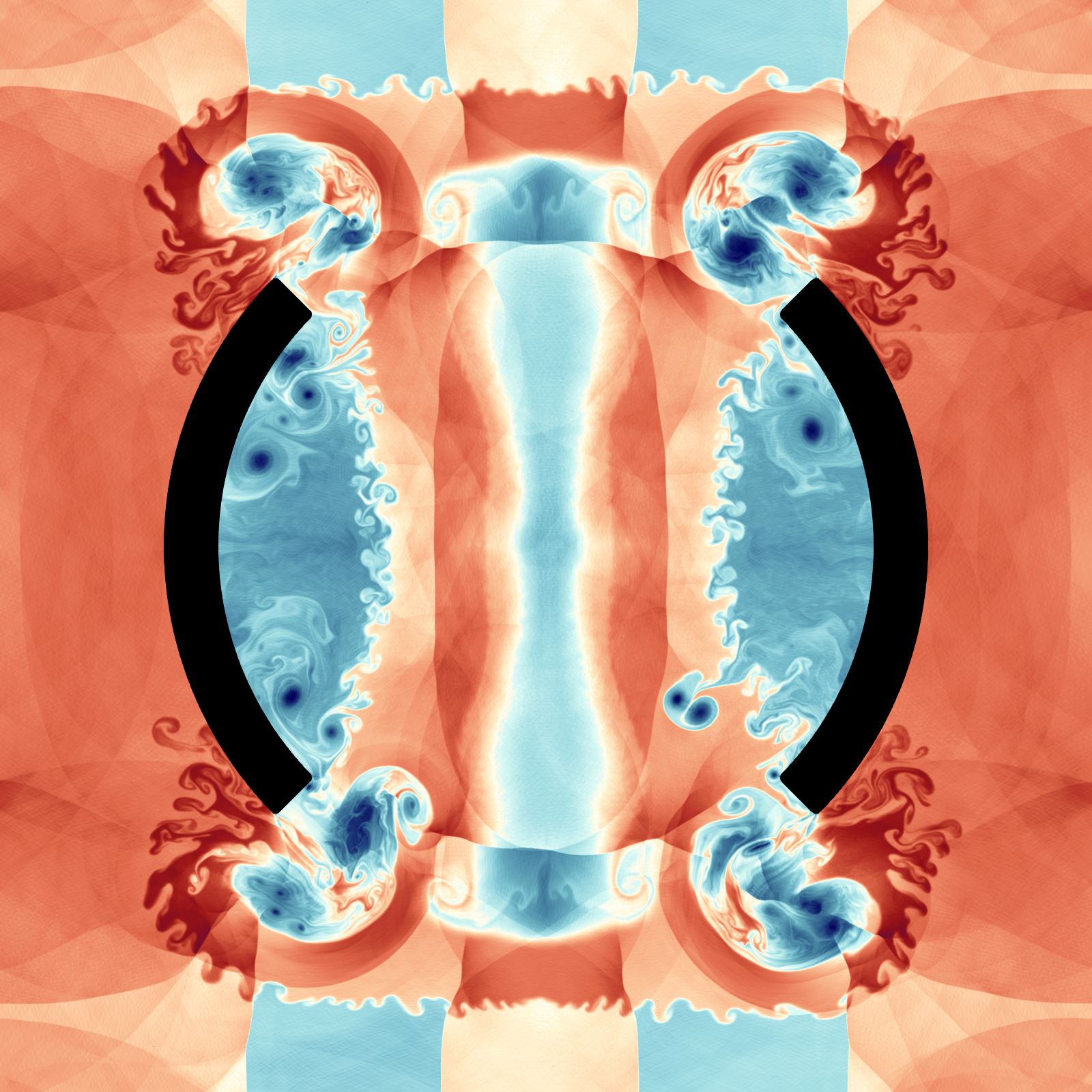}}}
        ~
        \subfloat[$t = 3$]{\label{fig:ant_2} \adjustbox{width=0.356\linewidth,valign=b}{\includegraphics[width=\textwidth]{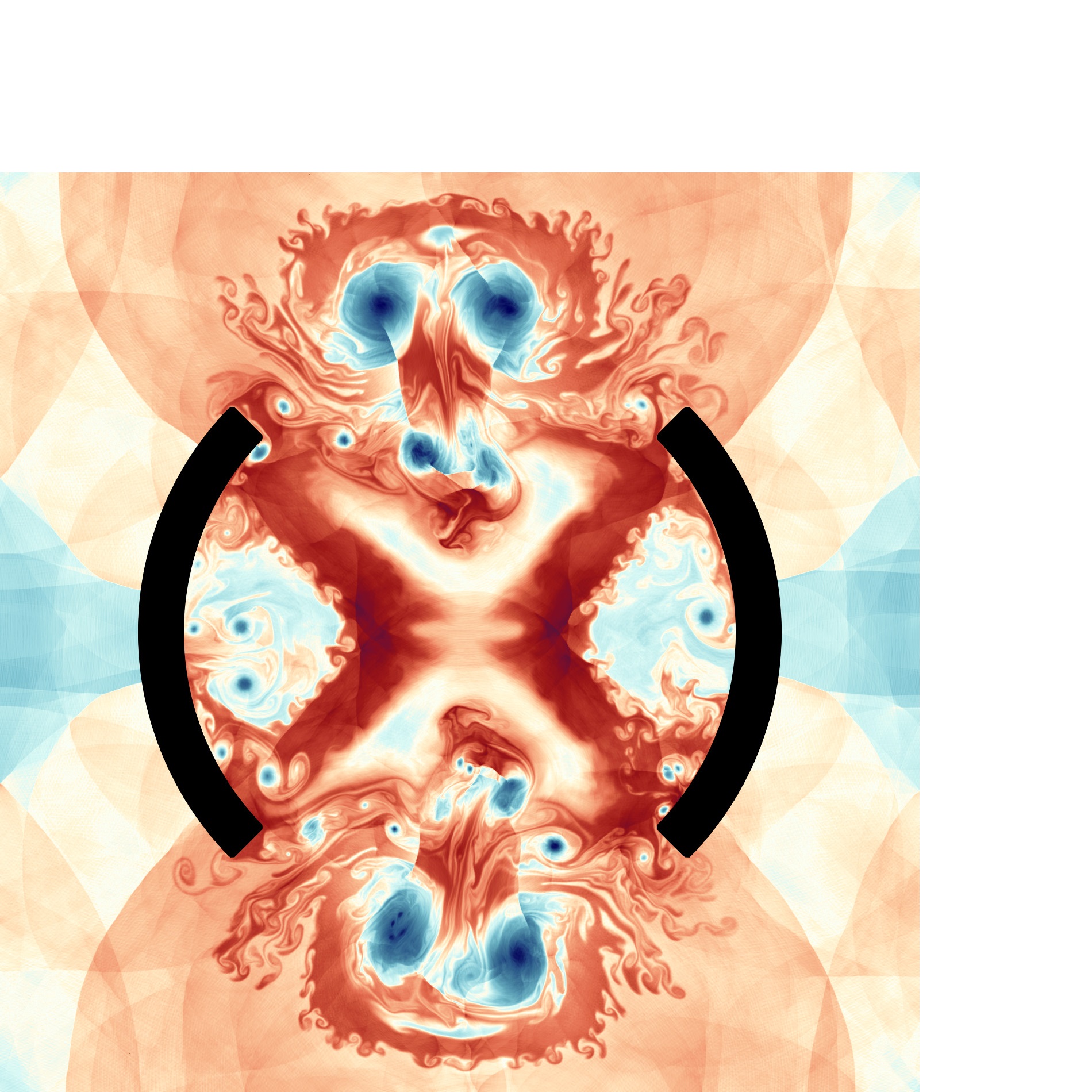}}}
        \hfill
        \newline
        \subfloat[$t = 1$]{\label{fig:antf_1} \adjustbox{width=0.3\linewidth,valign=b}{\includegraphics[width=\textwidth]{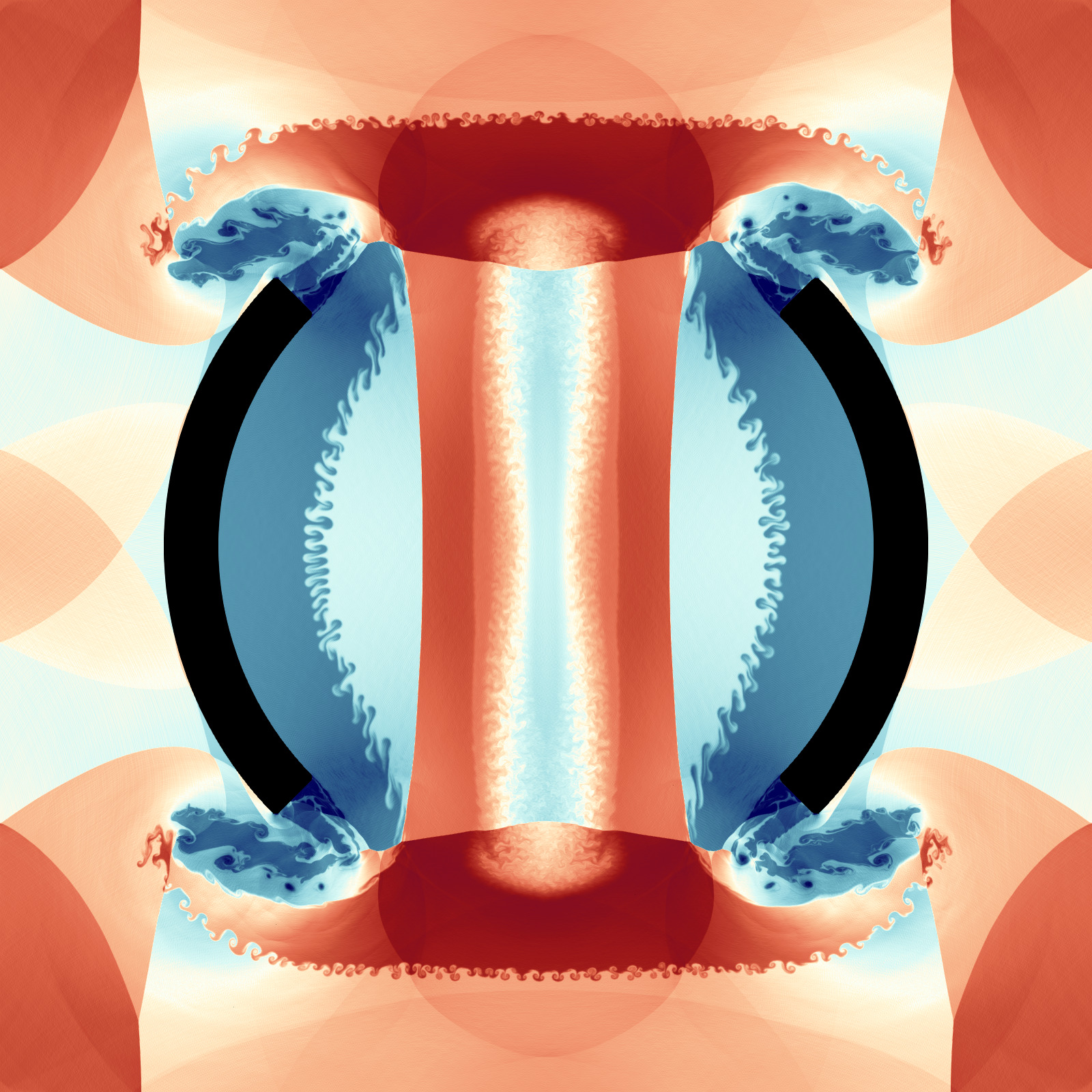}}}
        ~
        \subfloat[$t = 2$]{\label{fig:antf_2} \adjustbox{width=0.3\linewidth,valign=b}{\includegraphics[width=\textwidth]{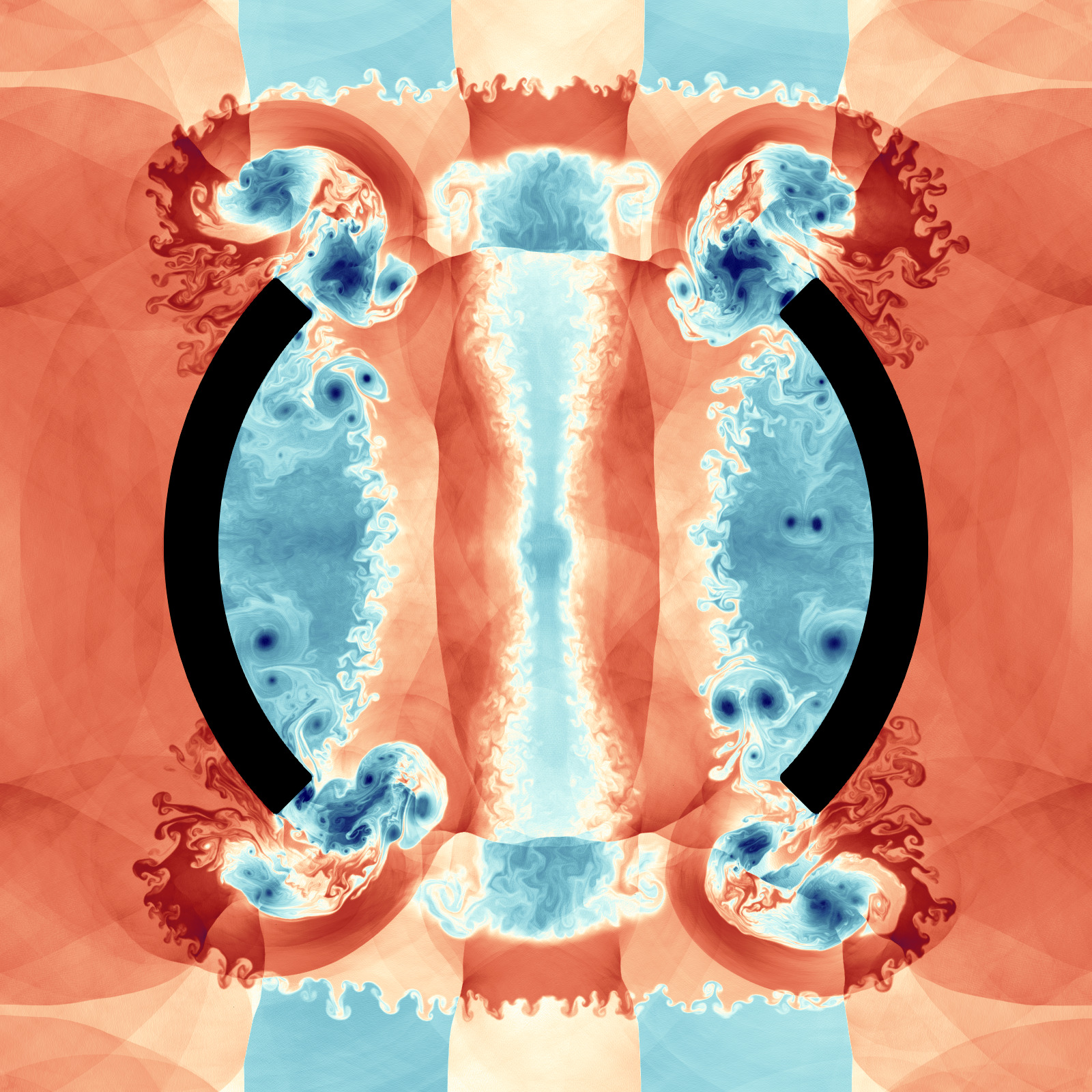}}}
        ~
        \subfloat[$t = 3$]{\label{fig:antf_3} \adjustbox{width=0.356\linewidth,valign=b}{\includegraphics[width=\textwidth]{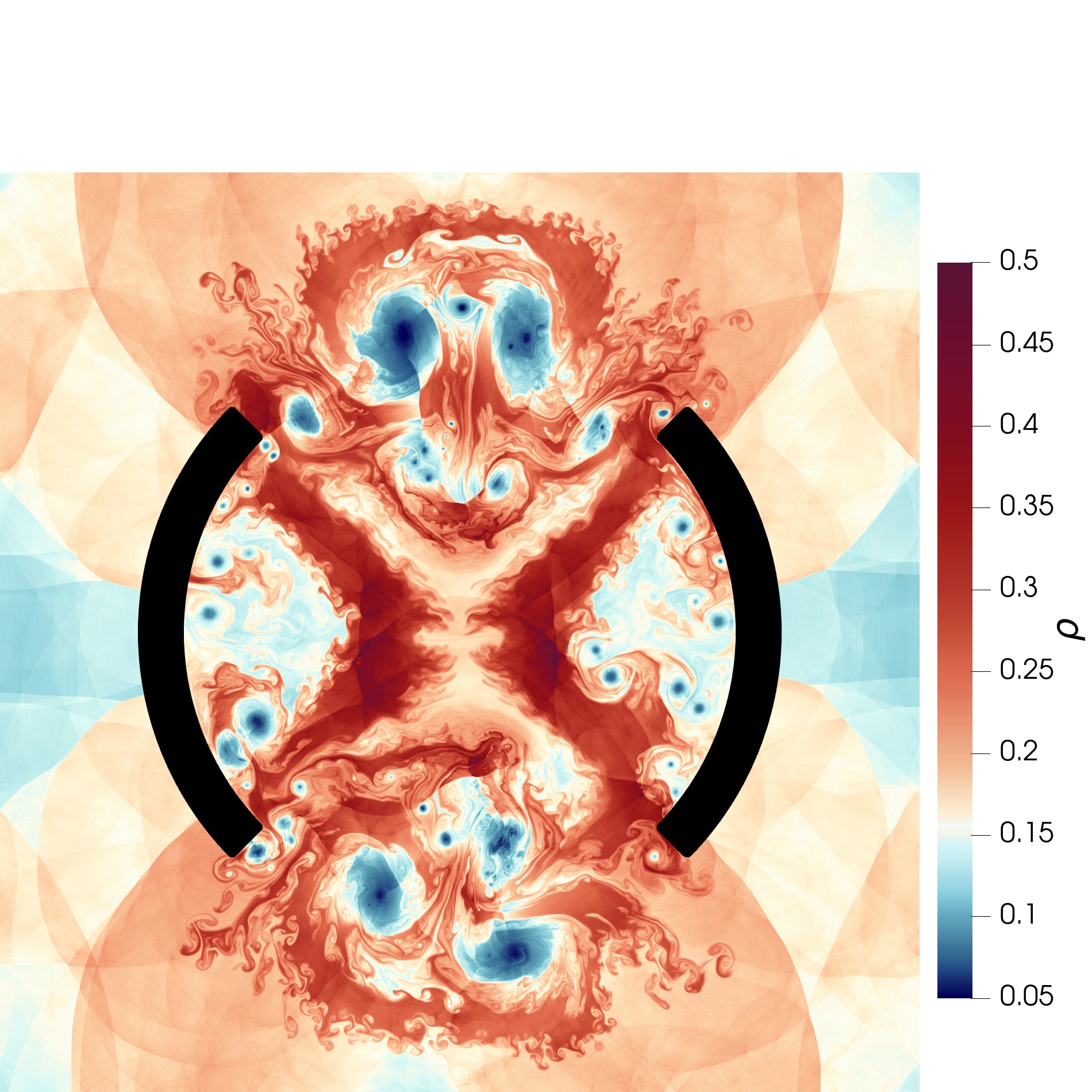}}}
        \hfill
        \newline
        \caption{\label{fig:pce} Logspaced contours of density for the partially-confined explosion problem at various time intervals computed using a $\mathbb P_4$ FR scheme with a coarse mesh (top row, $2\cdot10^5$ unstructured triangles) and a fine mesh (bottom row, $8\cdot10^5$ unstructured triangles).}
    \end{figure}

The results of a $\mathbb P_4$ FR approximation using the TR viscosity method ($c_\mu = 6$) are shown through the contours of density at various time intervals in \cref{fig:pce}. Initially, the shock front propagates radially outwards trailed by a slower-moving contact line while an expansion wave propagates towards the origin. The confinement region reflects the shock front, and the interaction of the reflected shock with the contact line seeds Richtmeyer-Meshkov-type instabilities which develop into complex vortical flow. With the TR viscosity method, the shock fronts were resolved with minimal spurious oscillations without excessively dissipating the vortical structures in the flow. Furthermore, the increase in mesh resolution resulted in a significantly better approximation of the fine-scale flow structures without additional tuning of the free parameter.

\subsubsection{ONERA M6}

\begin{figure}[htbp!]
    \centering
    \includegraphics[width=0.45\textwidth]{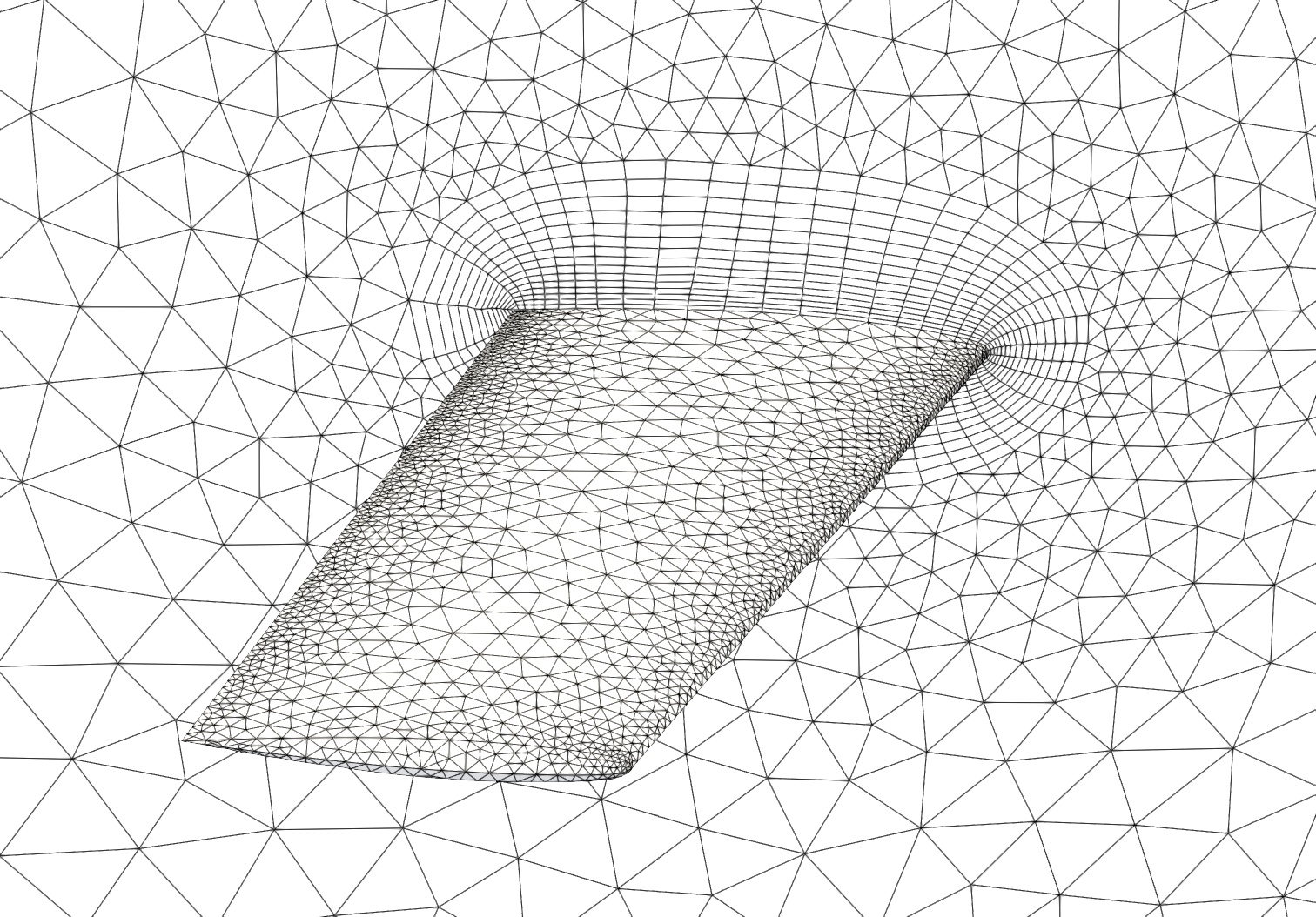}
    \caption{\label{fig:onera_mesh} Visualization of the mesh for the ONERA M6 wing. }
\end{figure}
For extensions to three-dimensional problems with applications in external aerodynamics, the inviscid transonic flow around the ONERA M6 wing was computed using the TR viscosity method. The geometry consists of a swept semi-span wing without twist utilizing a symmetric airfoil shape. The freestream Mach number is set to 0.8395 and the angle of attack is set to $3.06^{\circ}$. A second-order surface mesh of approximately $6\cdot10^3$ triangular elements was generated then extruded in the normal direction 16 times. The remaining mesh region consisted of unstructured tetrahedral elements. Slip-wall boundary conditions were applied to the wing surface and symmetry plane, and characteristic Riemann invariant boundary conditions were applied to the farfield. 

    \begin{figure}
        \centering
        \subfloat[Surface Pressure Coefficient]{\label{fig:onera_contours} \adjustbox{width=0.35\linewidth,valign=b}{\includegraphics[width=\textwidth]{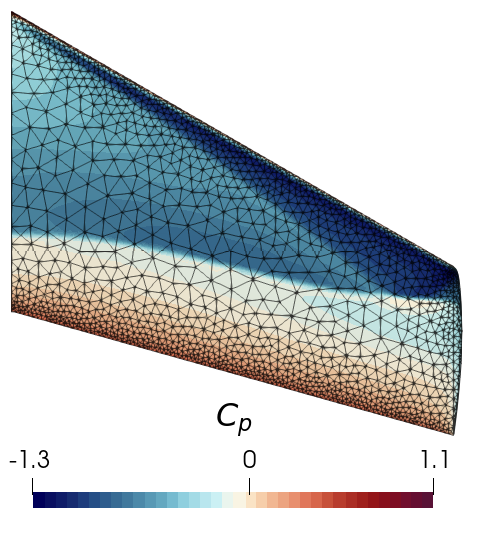}}}
        ~
        \subfloat[Pressure Coefficient Profile]{\label{fig:onera_cp} \adjustbox{width=0.45\linewidth,valign=b}{\input{./figs/onera_cpx0p65}}}
        \newline
        \caption{\label{fig:onera} Surface pressure coefficient contours with mesh overlay (left) and pressure coefficient profile at 65\% span (right) for the ONERA M6 wing computed using a $\mathbb P_4$ FR scheme. Experimental data of \citet{Schmitt1979} at a Reynolds number of $11{\cdot}10^6$ included for reference. }
    \end{figure}
    
The predictions of the pressure coefficient as computed by a $\mathbb P_4$ FR scheme with the TR viscosity approach ($c_\mu = 8$) are shown in \cref{fig:onera}. The surface pressure coefficient contours show the distinct "lambda-shock" pattern characteristic of this test case. Near the half-chord, the relatively coarse unstructured mesh was still able to sharply resolve the shock pattern, resolving it in the span of a few solution points. Similar observations were seen in the surface pressure coefficient profile at the cross-section correspond to 65\% of the span. The predicted profiles sharply resolved the pressure distribution without the introduction of spurious oscillations, showing satisfactory agreement with the \textit{viscous} experimental data of \citet{Schmitt1979}, performed at a Reynolds number of $11{\cdot}10^6$.

\subsection{Ideal Magnetohydrodynamics}
The TR viscosity method was then extended to the ideal magnetohydrodynamics (MHD) equations, written in the form of \cref{eq:hyp} as 
    \begin{equation}\label{eq:hyp_source}
        \partial_t \mathbf{u} (\mathbf{x}, t) + \boldsymbol{\nabla}{\cdot} \mathbf{F}\left(\mathbf{u}(\mathbf{x}, t)\right) = \mathbf{S}_{\boldsymbol{B}},
    \end{equation}
where the right-hand side includes additional source terms that are to be defined momentarily. Here, the solution and flux are represented by 
    \begin{equation}
        \mathbf{u} = \begin{bmatrix}
                \rho \\ \boldsymbol{\rho v} \\ \boldsymbol{B} \\ E
            \end{bmatrix}, \quad  \mathbf{F} = \begin{bmatrix}
                \boldsymbol{\rho v}\\
                \boldsymbol{\rho v}\otimes\mathbf{v} + \mathbf{I}\left(P + \frac{1}{2} \mathbf{B}{\cdot} \mathbf{B}\right) - \mathbf{B}\otimes\mathbf{B}\\
                \boldsymbol{v}\otimes\mathbf{B} - \mathbf{B}\otimes\boldsymbol{v}\\
            \left(E + P + \frac{1}{2} \mathbf{B}{\cdot} \mathbf{B}\right)\mathbf{v} - \mathbf{B}(\mathbf{v}{\cdot} \mathbf{B})
        \end{bmatrix},
    \end{equation}
where $\rho$ is the density, $\boldsymbol{\rho v}$ is the momentum, $E$ is the total energy, $P = (\gamma-1)\left(E - \shalf\rho\mathbf{v}{\cdot}\mathbf{v} - \shalf\mathbf{B}{\cdot}\mathbf{B}\right)$ is the pressure, $\mathbf{B}$ is the magnetic field, and $\gamma$ is the ratio of specific heat capacities. The symbol $\mathbf{I}$ denotes the identity matrix in $\mathbb{R}^{d\times d}$ and $\mathbf{v} = \boldsymbol{\rho v}/\rho$ denotes the velocity. For brevity, the solution is expressed in terms of a vector of primitive variables as $\mathbf{q}=[\rho,\mathbf{v},\mathbf{B},P]^T$. We additionally define a magnetic pressure as $P_b = \shalf(\gamma - 1)\mathbf{B}{\cdot}\mathbf{B}$. To enforce a solenoidal magnetic field, we utilize the approach of \citet{Powell1999} in which a source term proportional to the divergence of the magnetic field is included, written as
    \begin{equation}
        \mathbf{S}_{\boldsymbol{B}} = -\begin{bmatrix}
                0 \\ \boldsymbol{B} \\ \boldsymbol{u} \\ \boldsymbol{u}{\cdot}\boldsymbol{B}
            \end{bmatrix}\boldsymbol{\nabla}{\cdot}\boldsymbol{B},
    \end{equation}
where $\boldsymbol{\nabla}{\cdot}\boldsymbol{B}$ is calculated as the corrected divergence of the magnetic field using centered common interface values. 

\subsubsection{Rotor Problem}
The initial evaluation of the proposed method in the context of MHD was performed on the rotor problem of \citet{Balsara1999}, a two-dimensional test case for studying the effect of strong torsional Alfv\'{e}n waves relevant in applications such as star formation. The problem is defined on the domain $\Omega = [0,1]^2$ with the initial conditions 
    \begin{alignat*}{4}
        &\rho &&= 9\phi(r) + 1, &&\quad\quad  p && = 1,\\
        &u &&= -\frac{2}{r_0}(y - 0.5)\phi(r),  &&\quad\quad B_x && = \frac{5}{\sqrt{4 \pi}} ,\\
        &v &&= \frac{2}{r_0}(x - 0.5)\phi(r), &&\quad\quad  B_y && = 0 .
    \end{alignat*}
The blending between the inner and outer states was performed using a taper function $\phi(r)$, defined as 
    \begin{equation*}
        \phi(r) = \begin{cases}
            1, &\mbox{if } r \leq r_0, \\
            1 - \frac{r - r_0}{r_1 - r_0} &\mbox{if } r_0 < r \leq r_1,\\
            0, &\mbox{else},
        \end{cases} 
    \end{equation*}
where $r = \sqrt{x^2 + y^2}$, $r_0 = 0.1$, and $r_1 = 0.115$. A specific heat ratio of $\gamma = 1.4$ was used, and free boundary conditions were set on all sides of the domain. 
    
\begin{figure}[ht!]
        \centering
        \subfloat[Density]{ \adjustbox{width=0.3\linewidth,valign=b}{\includegraphics[width=\textwidth]{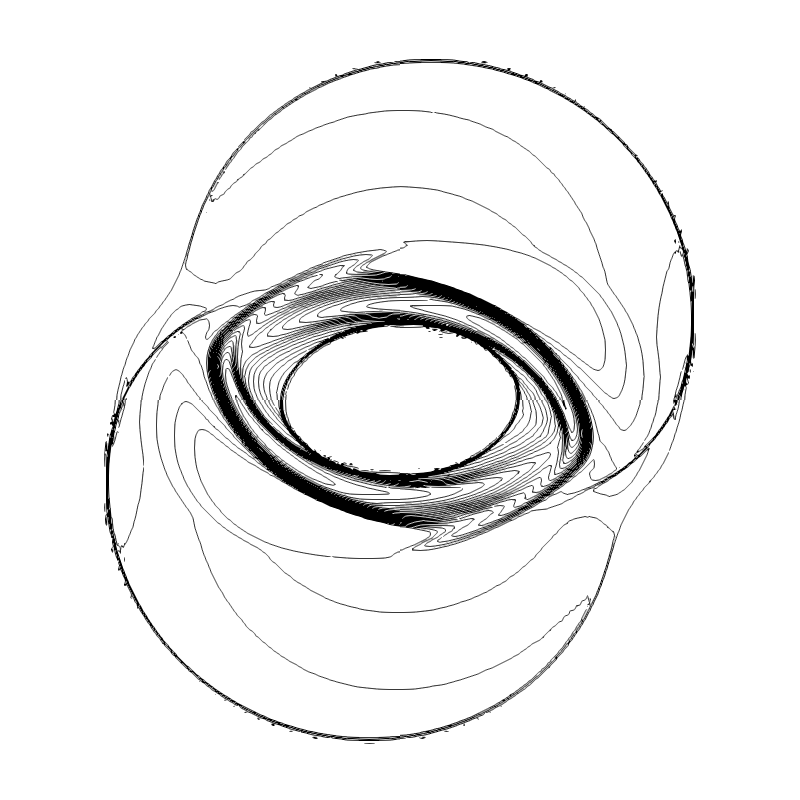}}}
        ~
        \subfloat[Pressure]{ \adjustbox{width=0.3\linewidth,valign=b}{\includegraphics[width=\textwidth]{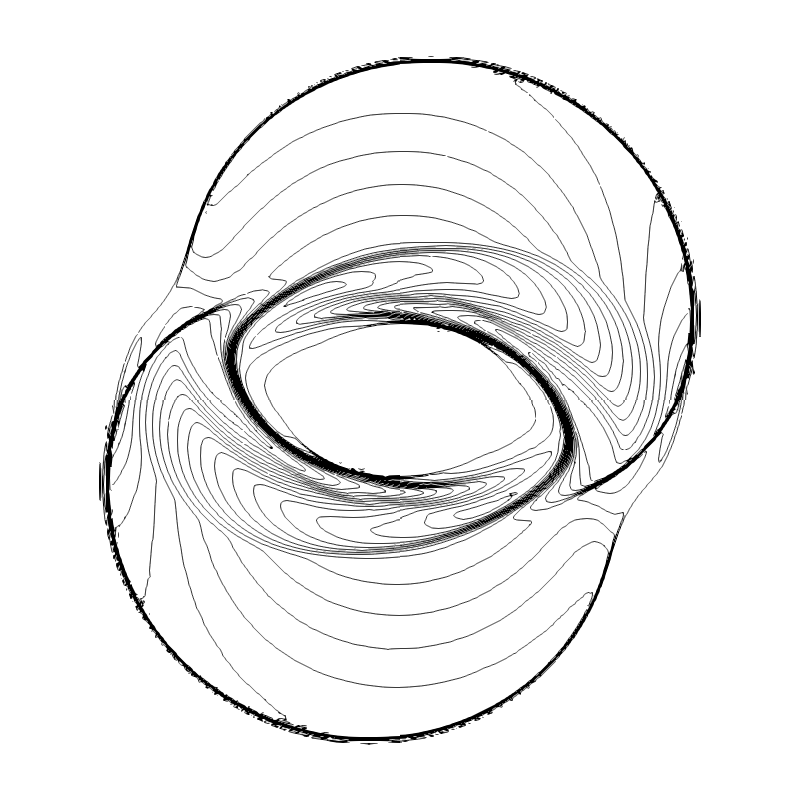}}}
        ~
        \subfloat[Magnetic Pressure]{ \adjustbox{width=0.3\linewidth,valign=b}{\includegraphics[width=\textwidth]{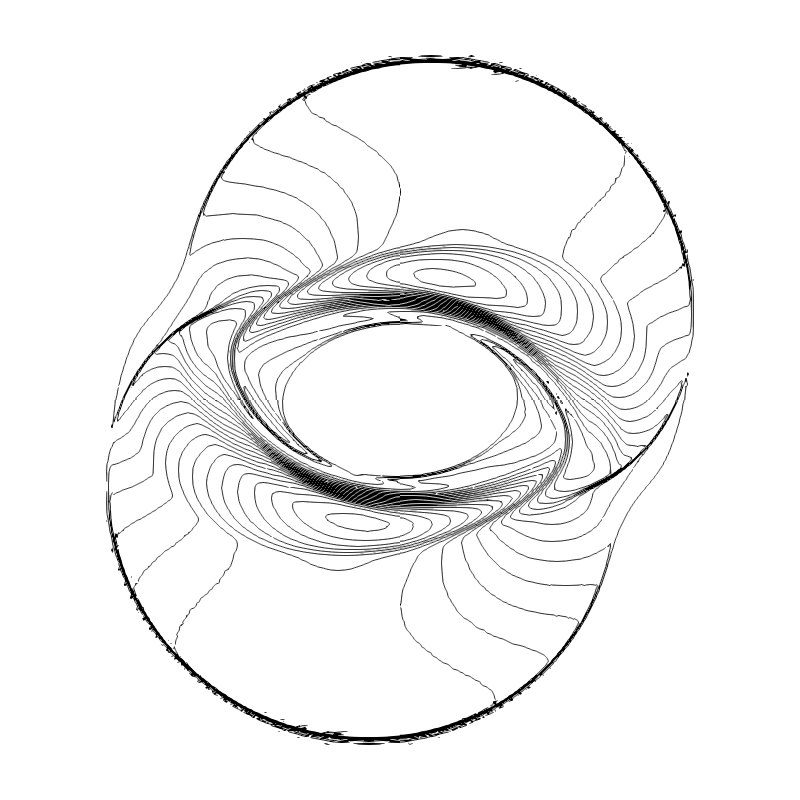}}}
        \caption{\label{fig:rotor2} Isocontours of density (left), pressure (middle), and magnetic pressure (right) for the rotor problem at $t = 0.15$ computed using a $\mathbb P_4$ FR scheme with $128^2$ elements.}
    \end{figure}

\begin{figure}[ht!]
        \centering
        \subfloat[Density]{\label{fig:rotor_d} \adjustbox{width=0.3\linewidth,valign=b}{\includegraphics[width=\textwidth]{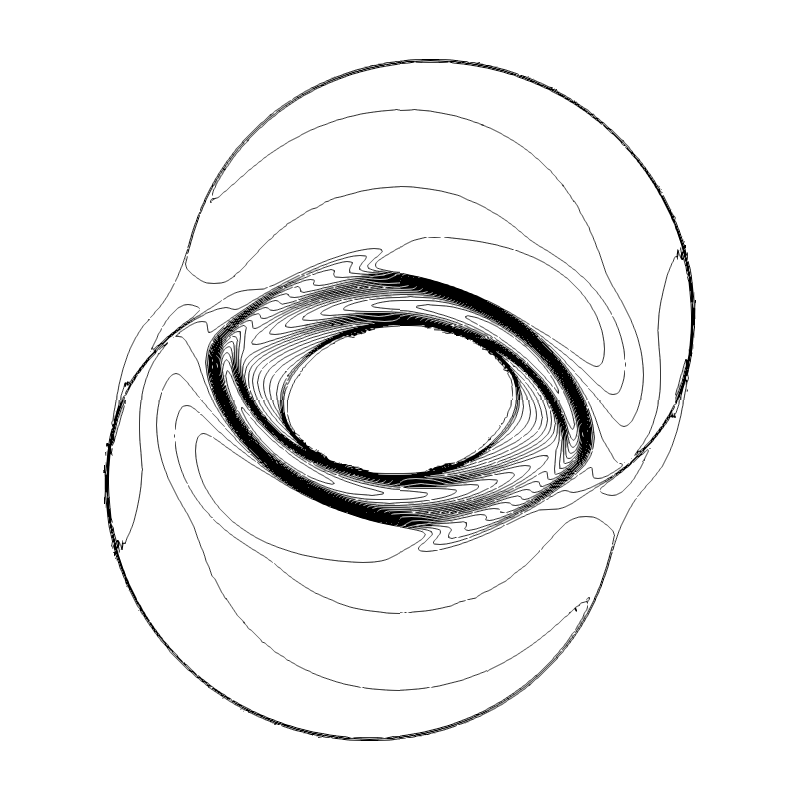}}}
        ~
        \subfloat[Pressure]{\label{fig:rotor_p} \adjustbox{width=0.3\linewidth,valign=b}{\includegraphics[width=\textwidth]{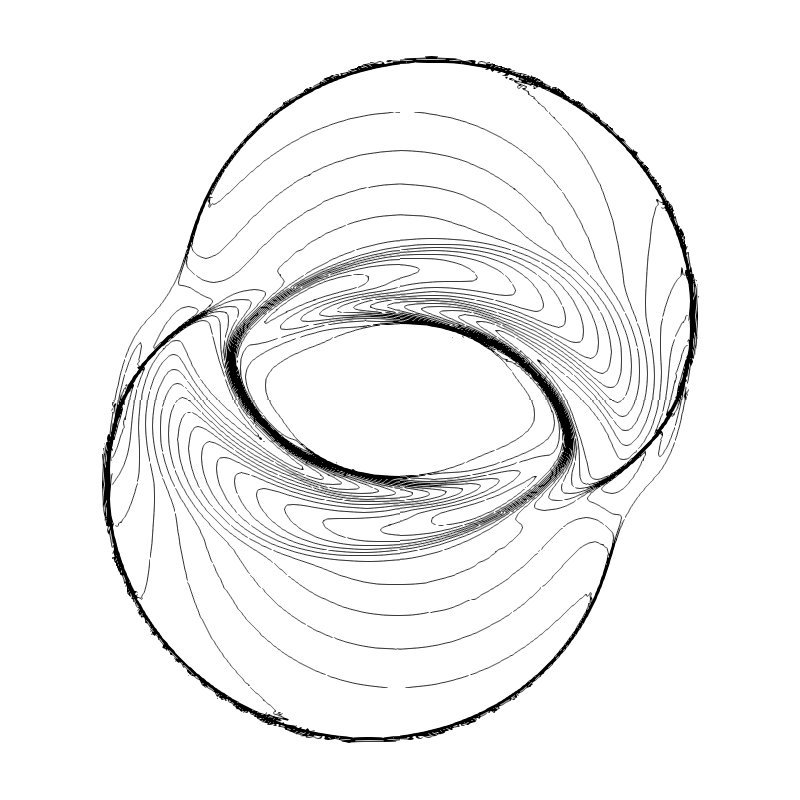}}}
        ~
        \subfloat[Magnetic Pressure]{\label{fig:rotor_mp} \adjustbox{width=0.3\linewidth,valign=b}{\includegraphics[width=\textwidth]{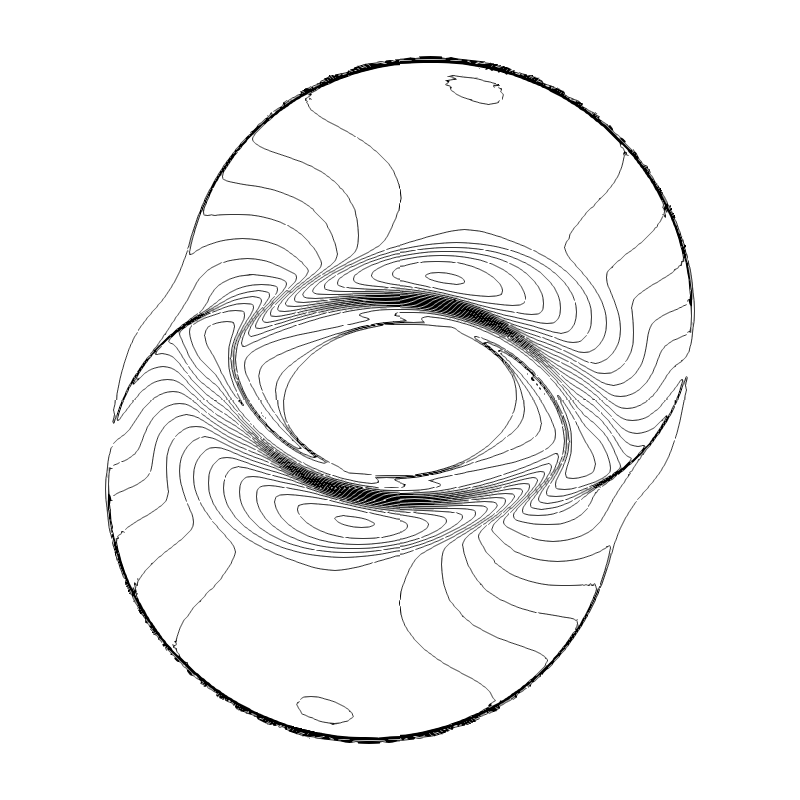}}}
        \caption{\label{fig:rotor} Isocontours of density (left), pressure (middle), and magnetic pressure (right) for the rotor problem at $t = 0.15$ computed using a $\mathbb P_9$ FR scheme with $64^2$ elements.}
    \end{figure}

The results of the TR viscosity method ($c_\mu = 2$) at $t = 0.15$ computed using a $\mathbb P_4$ FR scheme with $128^2$ quadrilateral elements and a $\mathbb P_9$ FR scheme with $64^2$ quadrilateral elements are shown in \cref{fig:rotor2} and \cref{fig:rotor}, respectively. The solution is shown through 20 logspaced isocontours of density on the range $[1, 10]$ as well as 20 equispaced isocontours of pressure and magnetic pressure on the ranges $[0, 2]$ and $[0, 1]$, respectively. For a fixed number of degrees of freedom, the isocontours from both the lower-order and higher-order approximations show excellent resolution of the flow and magnetic fields with effectively non-oscillatory profiles even in elements containing discontinuities. Furthermore, the discontinuities were resolved on a length scale on the order of 2-3 solution points, a small fraction of the total element size. No visibly noticeable deterioration of the results was observed when increasing the polynomial order of the approximation. 
    
    
\subsubsection{Orszag--Tang Vortex}
For a final test of the proposed approach, the efficacy of the TR viscosity method was evaluated on the canonical test case of \citet{Orszag1979}. The problem is defined on the periodic domain $\Omega = [0,1]^2$ with the initial conditions
    \begin{alignat*}{4}
        &\rho &&= \frac{25}{36} \pi, &&\quad\quad  p && = \frac{5}{12} \pi,\\
        &u &&= -\sin (2 \pi y),  &&\quad\quad B_x && = -\frac{1}{\sqrt{4 \pi}} \sin (2 \pi y),\\
        &v &&= \sin (2 \pi x), &&\quad\quad  B_y && = \frac{1}{\sqrt{4 \pi}} \sin (2 \pi x).
    \end{alignat*}
The specific heat ratio was set to $\gamma = 5/3$. The solution is initially smooth, but over time, the formation of MHD shocks and shock-shock interactions presents a challenge for high-resolution numerical schemes. The resulting flow fields are an ideal example of where vector-valued artificial viscosities offer benefits over their scalar-valued counterparts, as discontinuities may be present in the density field but not the magnetic field and vice versa.
    
    \begin{figure}[htbp!]
        \centering
        \subfloat[Density ($32^2$)]{\label{fig:OTV_d32} \adjustbox{width=0.4\linewidth,valign=b}{\includegraphics[width=\textwidth]{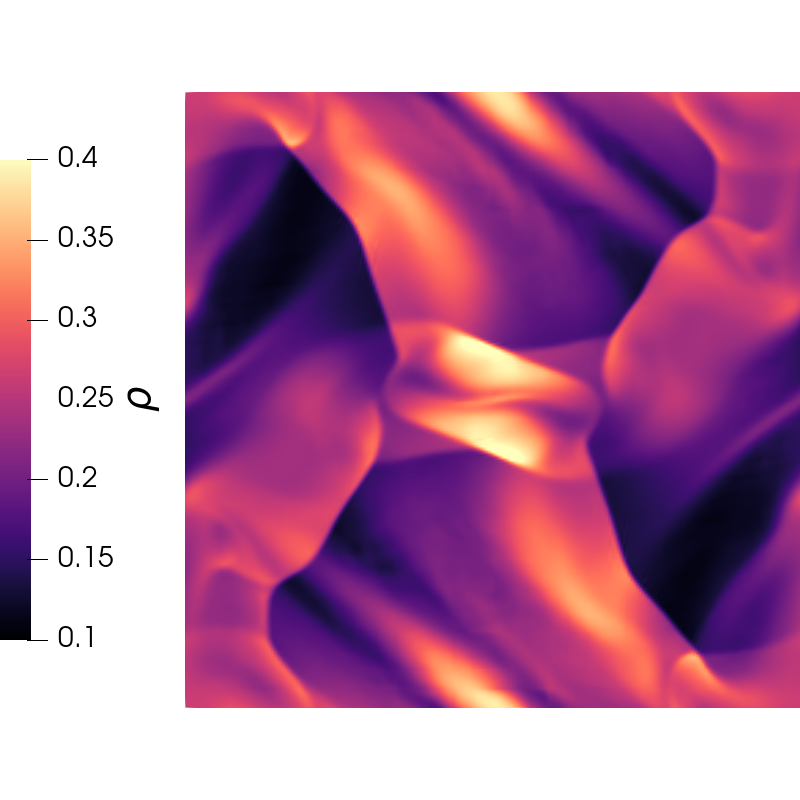}}}
        ~
        \subfloat[Density ($64^2$)]{\label{fig:OTV_d64} \adjustbox{width=0.4\linewidth,valign=b}{\includegraphics[width=\textwidth]{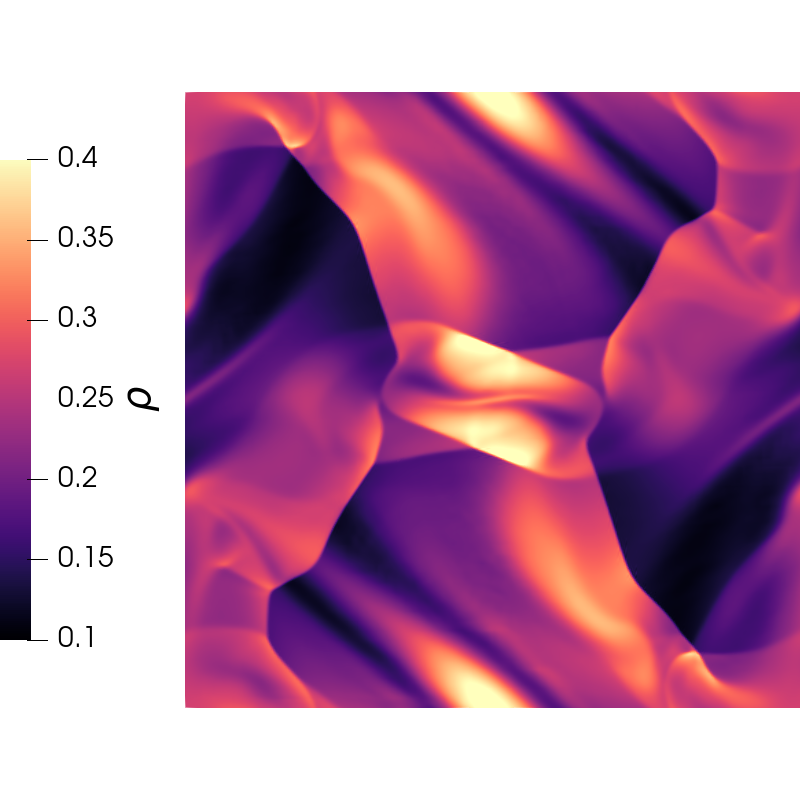}}}
        \newline
        \caption{\label{fig:OTV} Contours of density for the Orszag-Tang vortex at $t = 0.48$ computed using a $\mathbb P_7$ FR scheme with $32^2$ (left) and $64^2$ (right) elements.}
    \end{figure}
    
\begin{figure}
        \centering
        \subfloat[$\rho$-viscosity]{\label{fig:OTV_dv} \adjustbox{width=0.4\linewidth,valign=b}{\includegraphics[width=\textwidth]{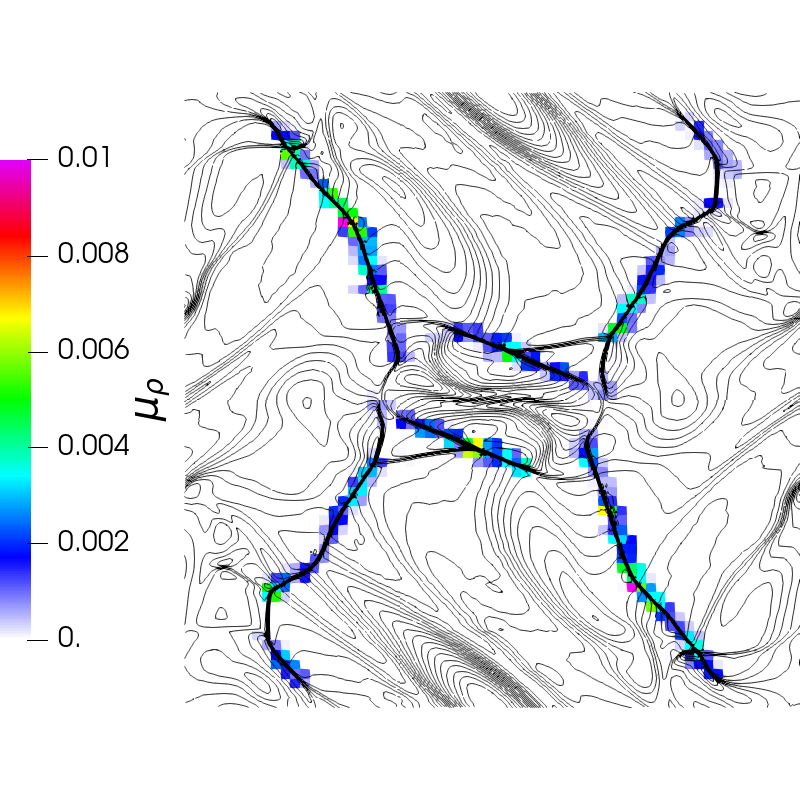}}}
        ~
        \subfloat[$B$-viscosity magnitude]{\label{fig:OTV_bv} \adjustbox{width=0.4\linewidth,valign=b}{\includegraphics[width=\textwidth]{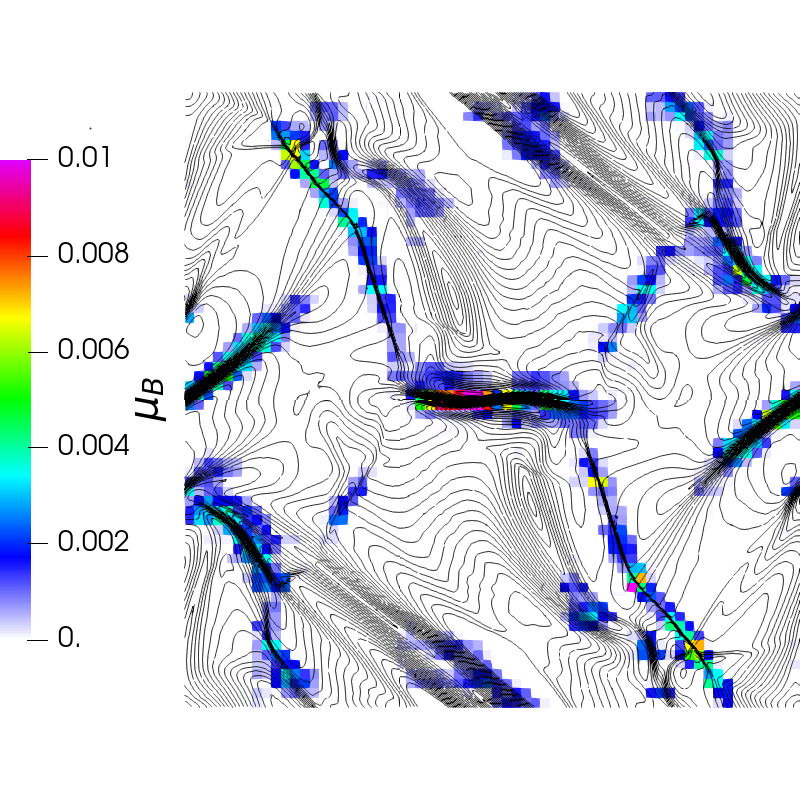}}}
        \newline
        \caption{\label{fig:OTV_v} Density component (left) and vector magnitude of the magnetic component (right) of the TR viscosity for the Orszag-Tang vortex at $t = 0.48$ computed using a $\mathbb P_7$ FR scheme with $64^2$ elements. Viscosity contours are overlayed with 20 equispaced isocontours of density and magnetic field magnitude, respectively.}
    \end{figure}

The contours of density at $t = 0.48$ computed using a $\mathbb P_7$ FR approximation with the TR viscosity method ($c_\mu = 3$) are shown in \cref{fig:pce} for meshes with $32^2$ and $64^2$ elements. The results show excellent resolution of the density field, even when utilizing a coarse mesh with a high approximation order. Without further tuning of the free parameter, the increase in mesh resolution resulted in significantly better prediction of the density field. To highlight the benefits of the vector-valued TR viscosity, the distributions of the density component and the vector magnitude of the magnetic components of the viscosity are shown in \cref{fig:OTV_v}. The viscosity distributions were particularly distinct between the two components, owing to the presence of hydrodynamic shocks and magnetic shocks independently of one another. As a result, artificial dissipation was only introduced where necessary for the various components. 
    
\begin{figure}[htbp!]
        \centering
        \adjustbox{width=0.7\linewidth,valign=b}{\input{./figs/otv_pressure}}
        \caption{\label{fig:OTV_p} Pressure profile of the Orszag-Tang vortex on the cross-section $y/L = 0.3125$ at $t = 0.48$ computed using a $\mathbb P_7$ FR scheme with $32^2$ and $64^2$ elements. Numerical results of \citet{Jiang1999} shown for reference.}
    \end{figure}
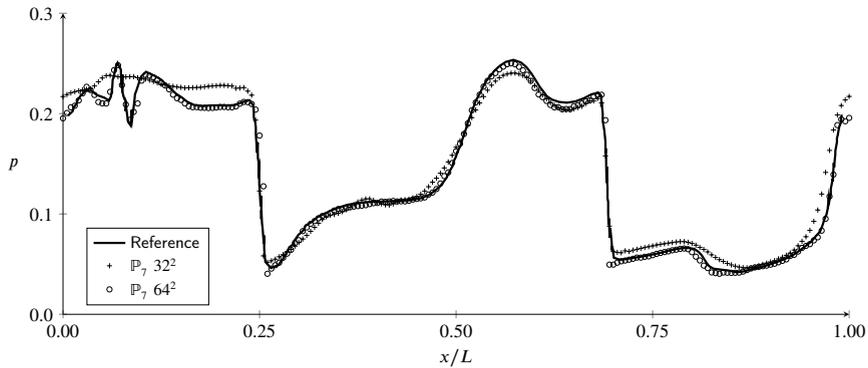

A comparison of the predicted pressure profile on the cross-section $y/L = 0.3125$ at $t = 0.48$ is shown in \cref{fig:OTV_p} with respect to the numerical results of \citet{Jiang1999}. On the coarse mesh, the discontinuities were still resolved with sub-element levels of resolution, but intricacies in the small-scale features of the pressure field were excessively diffused. When mesh resolution was increased, excellent agreement was observed between the predicted pressure field and the reference data, both with respect to the resolution of the discontinuities as well as the smooth regions and small-scale oscillations in the profile.  
    
%% file: figs/burgers_15.tex
    \begin{tikzpicture}[spy using outlines={rectangle, height=3cm,width=2.3cm, magnification=3, connect spies}]
		\begin{axis}[name=plot1,
		    axis line style={latex-latex},
		    axis x line=left,
            axis y line=left,
            clip mode=individual,
		    xlabel={$x$},
		    xtick={0,0.2,0.4,0.6,0.8,1},
    		xmin=0,
    		xmax=1,
    		x tick label style={
        		/pgf/number format/.cd,
            	fixed,
            	fixed zerofill,
            	precision=1,
        	    /tikz/.cd},
    		ylabel={$u$},
    		ylabel style={rotate=-90},
    		ytick={1, 1.5, 2.0, 2.5, 3.0},
    		ymin=1,
    		ymax=3,
    		y tick label style={
        		/pgf/number format/.cd,
            	fixed,
            	fixed zerofill,
            	precision=1,
        	    /tikz/.cd},
    		legend style={at={(0.97,.97)},anchor=north east,font=\small},
    		legend cell align={left},
    		style={font=\normalsize}]
    		
			\addplot[color=black, style={ thick}]
				table[x=x,y=u,col sep=comma,unbounded coords=jump]{./figs/data/burgers_exact.csv};
    		\addlegendentry{Exact}
    		
    		\addplot[color=black, style={ultra thin}, only marks, mark=o, mark options={scale=0.7}]
				table[x=x,y=u,col sep=comma,unbounded coords=jump]{./figs/data/burgers_p4_15elems.csv};
			\addlegendentry{$\mathbb P_{4}$}
			
            \addlegendimage{only marks, gray,mark=square*, fill=gray}
			\addlegendentry{$\mu$}
		\end{axis} 		

		\begin{axis}[name=plot2,
		    axis line style={latex-latex},
            axis y line=right,
            axis x line=none,
            clip mode=individual,
		    xtick={0,0.2,0.4,0.6,0.8,1},
    		xmin=0,
    		xmax=1,
    		x tick label style={
        		/pgf/number format/.cd,
            	fixed,
            	fixed zerofill,
            	precision=1,
        	    /tikz/.cd},
    		ylabel={$\mu$},
    		ylabel style={rotate=-90},
    		ytick={0, 4e-5, 8e-5, 12e-5, 16e-5, 2e-4},
    		ymin=0,
    		ymax=2e-4,
    		y tick label style={
        		/pgf/number format/.cd,
            	fixed,
            	fixed zerofill,
            	precision=1,
        	    /tikz/.cd},
    		style={font=\normalsize}]
    		
			\addplot[ybar, color=black, fill=black, fill opacity=0.1,
                draw opacity=0] 
				table[x=x,y=mu,col sep=comma,unbounded coords=jump]{./figs/data/burgers_p4_15elems.csv}; 
				
		\end{axis}

	\end{tikzpicture}

%% file: figs/burgers_5.tex
    \begin{tikzpicture}[spy using outlines={rectangle, height=3cm,width=2.3cm, magnification=3, connect spies}]
		\begin{axis}[name=plot1,
		    axis line style={latex-latex},
		    axis x line=left,
            axis y line=left,
            clip mode=individual,
		    xlabel={$x$},
		    xtick={0,0.2,0.4,0.6,0.8,1},
    		xmin=0,
    		xmax=1,
    		x tick label style={
        		/pgf/number format/.cd,
            	fixed,
            	fixed zerofill,
            	precision=1,
        	    /tikz/.cd},
    		ylabel={$u$},
    		ylabel style={rotate=-90},
    		ytick={1, 1.5, 2.0, 2.5, 3.0},
    		ymin=1,
    		ymax=3,
    		y tick label style={
        		/pgf/number format/.cd,
            	fixed,
            	fixed zerofill,
            	precision=1,
        	    /tikz/.cd},
    		legend style={at={(0.97,.97)},anchor=north east,font=\small},
    		legend cell align={left},
    		style={font=\normalsize}]
    		
			\addplot[color=black, style={ thick}]
				table[x=x,y=u,col sep=comma,unbounded coords=jump]{./figs/data/burgers_exact.csv};
    		\addlegendentry{Exact}
    		
    		\addplot[color=black, style={ultra thin}, only marks, mark=o, mark options={scale=0.7}]
				table[x=x,y=u,col sep=comma,unbounded coords=jump]{./figs/data/burgers_p14_5elems.csv};
			\addlegendentry{$\mathbb P_{14}$}
			
            \addlegendimage{only marks, gray,mark=square*, fill=gray}
			\addlegendentry{$\mu$}
		\end{axis} 		

		\begin{axis}[name=plot2,
		    axis line style={latex-latex},
            axis y line=right,
            axis x line=none,
            clip mode=individual,
		    xtick={0,0.2,0.4,0.6,0.8,1},
    		xmin=0,
    		xmax=1,
    		x tick label style={
        		/pgf/number format/.cd,
            	fixed,
            	fixed zerofill,
            	precision=1,
        	    /tikz/.cd},
    		ylabel={$\mu$},
    		ylabel style={rotate=-90},
    		ytick={0, 4e-5, 8e-5, 12e-5, 16e-5, 2e-4},
    		ymin=0,
    		ymax=2e-4,
    		y tick label style={
        		/pgf/number format/.cd,
            	fixed,
            	fixed zerofill,
            	precision=1,
        	    /tikz/.cd},
    		style={font=\normalsize}]
    		
			\addplot[ybar, color=black, fill=black, fill opacity=0.1,
                draw opacity=0] 
				table[x=x,y=mu,col sep=comma,unbounded coords=jump]{./figs/data/burgers_p14_5elems.csv};
				
		\end{axis}

	\end{tikzpicture}

%% file: figs/sod_density_p3.tex
    \begin{tikzpicture}[spy using outlines={rectangle, height=3cm,width=2.3cm, magnification=3, connect spies}]
		\begin{axis}[name=plot1,
		    axis line style={latex-latex},
		    axis x line=left,
            axis y line=left,
            clip mode=individual,
		    xlabel={$x$},
		    xtick={0,0.2,0.4,0.6,0.8,1},
    		xmin=0,
    		xmax=1,
    		x tick label style={
        		/pgf/number format/.cd,
            	fixed,
            	fixed zerofill,
            	precision=1,
        	    /tikz/.cd},
    		ylabel={$\rho$},
    		ylabel style={rotate=-90},
    		ytick={0,0.2,0.4,0.6,0.8,1},
    		ymin=0,
    		ymax=1.05,
    		y tick label style={
        		/pgf/number format/.cd,
            	fixed,
            	fixed zerofill,
            	precision=1,
        	    /tikz/.cd},
    		legend style={at={(0.03,0.03)},anchor=south west,font=\small},
    		legend cell align={left},
    		style={font=\normalsize}]
    		
			\addplot[color=black, style={thick}]
				table[x=x,y=r,col sep=comma,unbounded coords=jump]{./figs/data/sod_exact_reference.csv};
    		\addlegendentry{Exact}
    		
    		\addplot[color=black, style={ultra thin}, only marks, mark=o, mark options={scale=0.5}, mark repeat = 2, mark phase =0]
				table[x=x,y=rho,col sep=comma,unbounded coords=jump]{./figs/data/sod_p3_200dof_rho.csv};
			\addlegendentry{$\mathbb P_3$}
			
            \addlegendimage{only marks, gray,mark=square*, fill=gray}
			\addlegendentry{$\mu_\rho$}
		\end{axis} 		

		\begin{axis}[name=plot2,
		    axis line style={latex-latex},
            axis y line=right,
            axis x line=none,
            clip mode=individual,
		    xtick={0,0.2,0.4,0.6,0.8,1},
    		xmin=0,
    		xmax=1,
    		x tick label style={
        		/pgf/number format/.cd,
            	fixed,
            	fixed zerofill,
            	precision=1,
        	    /tikz/.cd},
    		ylabel={$\mu_{\rho}$},
    		ylabel style={rotate=-90},
    		ytick={0, 4e-5, 8e-5, 12e-5, 16e-5, 2e-4},
    		ymin=0,
    		ymax=2e-4,
    		y tick label style={
        		/pgf/number format/.cd,
            	fixed,
            	fixed zerofill,
            	precision=1,
        	    /tikz/.cd},
    		style={font=\normalsize}]
    		
			\addplot[ybar, color=black, fill=black, fill opacity=0.1,
                draw opacity=0] 
				table[x=x,y=mu_rho,col sep=comma,unbounded coords=jump]{./figs/data/sod_p3_200dof_murho.csv};
				
		\end{axis}

	\end{tikzpicture}

%% file: figs/sod_density_p9.tex
    \begin{tikzpicture}[spy using outlines={rectangle, height=3cm,width=2.3cm, magnification=3, connect spies}]
		\begin{axis}[name=plot1,
		    axis line style={latex-latex},
		    axis x line=left,
            axis y line=left,
		    xlabel={$x$},
		    xtick={0,0.2,0.4,0.6,0.8,1},
    		xmin=0,
    		xmax=1,
    		x tick label style={
        		/pgf/number format/.cd,
            	fixed,
            	fixed zerofill,
            	precision=1,
        	    /tikz/.cd},
    		ylabel={$\rho$},
    		ylabel style={rotate=-90},
    		ytick={0,0.2,0.4,0.6,0.8,1},
    		ymin=0,
    		ymax=1.05,
    		y tick label style={
        		/pgf/number format/.cd,
            	fixed,
            	fixed zerofill,
            	precision=1,
        	    /tikz/.cd},
    		legend style={at={(0.03,0.03)},anchor=south west,font=\small},
    		legend cell align={left},
    		style={font=\normalsize},
    		bar width=5pt]
    		
			\addplot[color=black, style={thick}]
				table[x=x,y=r,col sep=comma,unbounded coords=jump]{./figs/data/sod_exact_reference.csv};
    		\addlegendentry{Exact}
    		
    		\addplot[color=black, style={ultra thin}, only marks, mark=o, mark options={scale=0.5}, mark repeat = 2, mark phase =0]
				table[x=x,y=rho,col sep=comma,unbounded coords=jump]{./figs/data/sod_p9_200dof_rho.csv};
			\addlegendentry{$\mathbb P_9$}

            \addlegendimage{only marks, gray,mark=square*, fill=gray}
			\addlegendentry{$\mu_\rho$}

		\end{axis} 		
		\begin{axis}[name=plot2,
		    axis line style={latex-latex},
            axis y line=right,
            axis x line=none,
            clip mode=individual,
		    xtick={0,0.2,0.4,0.6,0.8,1},
    		xmin=0,
    		xmax=1,
    		x tick label style={
        		/pgf/number format/.cd,
            	fixed,
            	fixed zerofill,
            	precision=1,
        	    /tikz/.cd},
    		ylabel={$\mu_{\rho}$},
    		ylabel style={rotate=-90},
    		ytick={0, 4e-5, 8e-5, 12e-5, 16e-5, 2e-4},
    		ymin=0,
    		ymax=2e-4,
    		y tick label style={
        		/pgf/number format/.cd,
            	fixed,
            	fixed zerofill,
            	precision=1,
        	    /tikz/.cd},
    		style={font=\normalsize}]
    		
			\addplot[ybar, color=black, fill=black, fill opacity=0.1,
                draw opacity=0] 
				table[x=x,y=mu_rho,col sep=comma,unbounded coords=jump]{./figs/data/sod_p9_200dof_murho.csv};
				
		\end{axis}

	\end{tikzpicture}

%% file: figs/shu_density.tex
     \begin{tikzpicture}[spy using outlines={rectangle, height=3cm,width=2.3cm, magnification=3, connect spies}]
		\begin{axis}[name=plot1,
		    axis line style={latex-latex},
		    axis x line=left,
            axis y line=left,
		    xlabel={$x$},
		    xtick={-5,-2.5,0,2.5,5},
    		xmin=-5,
    		xmax=5,
    		x tick label style={
        		/pgf/number format/.cd,
            	fixed,
            	fixed zerofill,
            	precision=1,
        	    /tikz/.cd},
    		ylabel={$\rho$},
    		ylabel style={rotate=-90},
    		ytick={0,1,2,3,4,5},
    		ymin=0,
    		ymax=5,
    		y tick label style={
        		/pgf/number format/.cd,
            	fixed,
            	fixed zerofill,
            	precision=0,
        	    /tikz/.cd},
    		legend style={at={(0.03,0.075)},anchor=south west,font=\small},
    		legend cell align={left},
    		style={font=\normalsize}]
    		
			\addplot[color=black, style={thin}]
				table[x=x,y=d,col sep=comma,unbounded coords=jump]{./figs/data/osher_p0_2000.csv};
    		\addlegendentry{Reference}

			\addplot[color=black, style={ultra thin}, only marks, mark=o, mark options={scale=0.5}, mark repeat = 2, mark phase =0]
				table[x expr={\thisrow{x}-5},y=rho,col sep=comma,unbounded coords=jump]{./figs/data/shuosher_p4_500dof_roe_cmu5.csv};
    		\addlegendentry{$\mathbb P_4$}
    		
            \addlegendimage{only marks, gray,mark=square*, fill=gray}
			\addlegendentry{$\mu_\rho$}

		\end{axis} 		
		\begin{axis}[name=plot2,
		    axis line style={latex-latex},
            axis y line=right,
            axis x line=none,
            clip mode=individual,
		    xtick={},
    		xmin=-5,
    		xmax=5,
    		x tick label style={
        		/pgf/number format/.cd,
            	fixed,
            	fixed zerofill,
            	precision=1,
        	    /tikz/.cd},
    		ylabel={$\mu_{\rho}$},
    		ylabel style={rotate=-90},
    		ytick={0, 2e-4, 4e-4, 6e-4, 8e-4, 1e-3},
    		ymin=0,
    		ymax=1e-3,
    		y tick label style={
        		/pgf/number format/.cd,
            	fixed,
            	fixed zerofill,
            	precision=1,
        	    /tikz/.cd},
    		style={font=\normalsize}]
    		
			\addplot[ybar, color=black, fill=black, fill opacity=0.1,
                draw opacity=0] 
				table[x expr={\thisrow{x}-5},y=mu_rho,col sep=comma,unbounded coords=jump]{./figs/data/shuosher_viscosity.csv};
				
		\end{axis} 	
		
	\end{tikzpicture}

%% file: figs/exp_geo.tex
     \begin{tikzpicture}[spy using outlines={rectangle, height=3cm,width=2.3cm, magnification=3, connect spies}]
		\begin{axis}[name=plot1,
		    axis line style={draw=none},
		    tick style={draw=none},
		    axis x line=left,
            axis y line=left,
            axis equal image,
            clip mode=individual,
    		xmin=-1.1,
    		xmax=1.1,
    		xticklabels={,,},
    		ymin=-1.1,
    		ymax=1.1,
    		yticklabels={,,},
    		style={font=\small},
    		scale = 1]

        \draw[-, red, fill=red!10] (axis cs:0,0) circle[black, radius=0.3];
        
		\draw[-] (axis cs:-1.0, -1.0) -- (axis cs:1.0, -1.0);
		\draw[-] (axis cs:1.0, -1.0) -- (axis cs:1.0, 1.0);
		\draw[-] (axis cs:1.0, 1.0) -- (axis cs:-1.0, 1.0);
		\draw[-] (axis cs:-1.0, 1.0) -- (axis cs:-1.0, -1.0);

        \draw[name path=A] (axis cs:0.6, 0.0) arc (0:45:.6);
        \draw[name path=B] (axis cs:0.7, 0.0) arc (0:45:.7);
        \draw[name path=C] (axis cs:0.6, 0.0) arc (0:-45:.6);
        \draw[name path=D] (axis cs:0.7, 0.0) arc (0:-45:.7);
		\draw[-] (axis cs:0.424264069, 0.424264069) -- (axis cs:0.494974747, 0.494974747);
		\draw[-] (axis cs:0.424264069, -0.424264069) -- (axis cs:0.494974747, -0.494974747);
		\tikzfillbetween[of=A and B]{gray!40};
		\tikzfillbetween[of=C and D]{gray!40};
        
        \draw[name path=E] (axis cs:-0.6, 0.0) arc (180:225:.6);
        \draw[name path=F] (axis cs:-0.7, 0.0) arc (180:225:.7);
        \draw[name path=G] (axis cs:-0.6, 0.0) arc (180:135:.6);
        \draw[name path=H] (axis cs:-0.7, 0.0) arc (180:135:.7);
		\draw[-] (axis cs:-0.424264069, 0.424264069) -- (axis cs:-0.494974747, 0.494974747);
		\draw[-] (axis cs:-0.424264069, -0.424264069) -- (axis cs:-0.494974747, -0.494974747);
		\tikzfillbetween[of=E and F]{gray!40};
		\tikzfillbetween[of=G and H]{gray!40};

	    \node at (axis cs:0.0,0.15) 
	    {\footnotesize $\mathbf{u}_L$};
	    
	    \node at (axis cs:0.0,-0.75) 
	    {\footnotesize $\mathbf{u}_R$};
	    
		\draw[->, densely dotted] (axis cs:0, 0) -- (axis cs:0, -0.3);
	    \node at (axis cs:-.05, -0.15) 
	    {\tiny $r_1$};

		\draw[->, densely dotted] (axis cs:0, 0) -- (axis cs:0.55432772, 0.229610059);
	    \node at (axis cs:0.4, 0.23) 
	    {\tiny $r_2$};
		\draw[->, densely dotted] (axis cs:0, 0) -- (axis cs:0.646715673, -0.267878402);
	    \node at (axis cs:0.4, -0.23) 
	    {\tiny $r_3$};
		
		\draw[-, densely dotted] (axis cs:0, 0) -- (axis cs:-0.424264069, 0.424264069);
		\draw[-, densely dotted] (axis cs:0, 0) -- (axis cs:-0.424264069, -0.424264069);
        \draw[<->, densely dotted] (axis cs:-0.318198051, -0.318198051) arc (225:135:.45);
	    \node at (axis cs:-0.5, 0.0) 
	    {\tiny $\theta$};

		\draw[->, densely dotted] (axis cs:0.55, 0.60) -- (axis cs:0.494974747, 0.494974747);
	    \node at (axis cs:0.55, 0.65) 
	    {\tiny $r_f \ (\times 8)$};

		\end{axis}

	\end{tikzpicture}

%% file: figs/onera_cpx0p65.tex
    \begin{tikzpicture}[spy using outlines={rectangle, height=3cm,width=2.3cm, magnification=3, connect spies}]
		\begin{axis}[name=plot1,
		    axis line style={latex-latex},
		    axis x line=left,
            axis y line=left,
            clip mode=individual,
		    xlabel={$x/c$},
		    xtick={0,0.2,0.4,0.6,0.8,1},
    		xmin=-0.05,
    		xmax=1.05,
    	    x tick label style={
        		/pgf/number format/.cd,
            	fixed,
            	fixed zerofill,
            	precision=1,
        	    /tikz/.cd},
    		ylabel={$C_p$},
    		ylabel style={rotate=-90},
    		ytick={-1.5, -1.0, -0.5, 0.0, 0.5, 1.0},
    		ymin=-1.5,
    		ymax=1,
    		y dir=reverse,
    		y tick label style={
        		/pgf/number format/.cd,
            	fixed,
            	fixed zerofill,
            	precision=1,
        	    /tikz/.cd},
    		legend style={at={(1.05, 1.0)},anchor=north east,font=\small},
    		legend cell align={left},
    		style={font=\normalsize}]
    		
			\addplot[color=black, style={thin}, only marks, mark=*, mark options={scale=0.8}]
				table[x=x,y=cp,col sep=comma,unbounded coords=jump]{./figs/data/onera_refdata_x0p65.csv};
    		\addlegendentry{Experiment}
    		
			\addplot[color=black, style={thick}]
				table[x=x,y=cp,col sep=comma,unbounded coords=jump]{./figs/data/onera_cp_x0p65.csv};
    		\addlegendentry{$\mathbb P_4$}

		\end{axis} 		
	\end{tikzpicture}

%% file: figs/otv_pressure.tex
    \begin{tikzpicture}[spy using outlines={rectangle, height=3cm,width=2.3cm, magnification=3, connect spies}]
		\begin{axis}[name=plot1,
            width=\textwidth,
            height=\axisdefaultheight,
		    axis line style={latex-latex},
		    axis x line=left,
            axis y line=left,
            clip mode=individual,
		    xlabel={$x/L$},
		    xtick={0, 0.25, 0.5, 0.75, 1},
    		xmin=0,
    		xmax=1,
    	    x tick label style={
        		/pgf/number format/.cd,
            	fixed,
            	fixed zerofill,
            	precision=2,
        	    /tikz/.cd},
    		ylabel={$p$},
    		ylabel style={rotate=-90},
    		ytick={0, 0.10, 0.20,0.30},
    		ymin=0,
    		ymax=0.3,
    		y tick label style={
        		/pgf/number format/.cd,
            	fixed,
            	fixed zerofill,
            	precision=1,
        	    /tikz/.cd},
    		legend style={at={(0.03,0.03)},anchor=south west,font=\small},
    		legend cell align={left},
    		style={font=\normalsize}]
    		

			\addplot[color=black, style={very thick}]
				table[x index=0,y index=1,col sep=comma,unbounded coords=jump]{./figs/data/ref_data_jiang_wu_1999.csv};
			\addlegendentry{ Reference}
				
			\addplot[color=black, style={ultra thin}, only marks, mark=+, mark options={scale=0.7}]
				table[x=Points_0,y=Pressure,col sep=comma,unbounded coords=jump]{./figs/data/p7_32x32_y0p3125.csv};
			\addlegendentry{$\mathbb P_7$ $32^2$}
				
			\addplot[color=black, style={ultra thin}, only marks, mark=o, mark options={scale=0.7}]
				table[x=Points_0,y=Pressure,col sep=comma,unbounded coords=jump]{./figs/data/p7_64x64_y0p3125.csv};
			\addlegendentry{$\mathbb P_7$ $64^2$}

		\end{axis} 		
	\end{tikzpicture}

%% file: conclusions.tex
\section{Conclusion}\label{sec:conclusion}
A novel approach for formulating an artificial viscosity for shock capturing in nonlinear hyperbolic conservation laws is introduced in this work. The proposed method relies on the notion of time-reversibility in the context of hyperbolic conservation laws and utilizes the property that hyperbolic systems are irreversible in the vicinity of a shock. By calculating the ability of the scheme to recover a given state under reverse temporal integration, a nonlinear stabilization technique is formed by injecting a proportionate amount of artificial diffusion in the regions where the system is (numerically) irreversible. The calculation of this artificial viscosity does not require any additional governing equations or \textit{a priori} knowledge of the system, is independent of the mesh and approximation order, and results in one tunable parameter. The primary novelty of the present work is that for systems of equations, the resulting artificial viscosity is unique for each component of the system which can be especially beneficial for systems in which some components exhibit discontinuities while others do not.
To assess the efficacy of the proposed approach, a series of multi-dimensional nonlinear hyperbolic problems were evaluated using a high-order discontinuous spectral element method on unstructured grids.
Even with very high order approximations, sub-element resolution of discontinuities was generally observed without the introduction of excessive diffusion in smooth regions. These findings indicate that the proposed approach is a simple yet effective and robust shock capturing method for arbitrary nonlinear systems of hyperbolic conservation laws. Extensions of the methods in this work will focus on applications to mixed hyperbolic-parabolic systems such as the Navier--Stokes equations and alternate non-parabolic forms of viscous regularization.